\documentclass[a4paper,oneside,11pt]{article}
 
\usepackage[a4paper]{geometry}
\usepackage{aeguill}                 
\usepackage{graphicx}
\usepackage{caption}
\usepackage{amsmath,amsfonts,amssymb,amsthm}
\usepackage{mathabx}
\usepackage{pstricks,comment} 
\usepackage{pst-plot}
\usepackage{pstricks-add}
\usepackage{multido}
\usepackage{url}
\usepackage{epstopdf,enumerate}
\usepackage{srcltx}
\usepackage[utf8]{inputenc} 
\usepackage[T1]{fontenc} 
\usepackage[english]{babel} 
\usepackage{hyperref}
\usepackage[all]{xy} 
\usepackage{titlesec}
\usepackage{mathabx}
\usepackage{tikz}
\usepackage{fancyhdr}
\usepackage[us,12hr]{datetime}
\usepackage{mathrsfs}
\fancypagestyle{plain}{
\fancyhf{}
\rfoot{ {\ddmmyyyydate\today}}
\lfoot{\thepage}
}

    \newtheoremstyle{TheoremNum}
        {\topsep}{\topsep}              %%% space between body and thm
        {\itshape}                      %%% Thm body font
        {}                              %%% Indent amount (empty = no indent)
        {\bfseries}                     %%% Thm head font
        {.}                             %%% Punctuation after thm head
        { }                             %%% Space after thm head
        {\thmname{#1}\thmnote{ \bfseries #3}}%%% Thm head spec

{\theoremstyle {definition} \newtheorem {defi} {Definition}[section]}
{\theoremstyle {plain}  \newtheorem {theo} [defi] {Theorem}}
{\theoremstyle {plain}  }
{\theoremstyle {plain} \newtheorem {prop} [defi] {Proposition}}
{\theoremstyle {plain} \newtheorem {lem}[defi] {Lemma}}
{\theoremstyle {plain} \newtheorem {rmq}[defi] {Remark}}
{\theoremstyle {plain} }
{\theoremstyle{TheoremNum} }
{\theoremstyle{TheoremNum} }
{\theoremstyle{TheoremNum} }

\newcommand{\Aut}{\mathrm{Aut}}
\newcommand{\Out}{\mathrm{Out}}

\newcommand{\Mod}{\mathrm{Mod}}

\newcommand{\Stab}{\mathrm{Stab}}

\newcommand{\Curr}{\mathrm{Curr}}
\newcommand{\PCurr}{\mathbb{P}\mathrm{Curr}}

\newcommand{\NN}{\mathbb{N}}
\newcommand{\RR}{\mathbb{R}}

\newcommand{\dem}{\noindent{\bf Proof. }}

\title{Currents relative to a malnormal subgroup system}
\author{Yassine Guerch}
\date{\today}
\begin{document}
\maketitle
\renewcommand*\labelenumi{(\theenumi)}

\begin{abstract}
This paper introduces a new topological space associated with a nonabelian free group $F_n$ of rank $n$ and a malnormal subgroup system $\mathcal{A}$ of $F_n$, called the space of currents relative to $\mathcal{A}$, which are $F_n$-invariant measures on an appropriate subspace of the double boundary of $F_n$. The extension from free factor systems as considered by Gupta to malnormal subgroup systems is necessary in order to fully study the growth under iteration of outer automorphisms of $F_n$, and requires the introduction of new techniques on cylinders. We in particular prove that currents associated with elements of $F_n$ which are not contained in a conjugate of a subgroup of $\mathcal{A}$ are dense in the space of currents relative to $\mathcal{A}$.
\end{abstract}

\section{Introduction}

Let $n \geq 2$. This paper is the first of a sequence of papers where we study the exponential growth of elements of $\Out(F_n)$, the outer automorphism group of a nonabelian free group $F_n$ of rank $n$. Let $[g]$ be the conjugacy class of a nontrivial element $g$ of $F_n$, let \mbox{$\phi \in \Out(F_n)$} and let $\Phi \in \Aut(F_n)$ be a representative of $\phi$. We say that $[g]$ has \emph{exponential growth under iterates of $\phi$} if there exists a basis $\mathcal{B}$ of $F_n$ such that the length of $[\Phi^n(g)]$ with respect to the word metric relative to $\mathcal{B}$ grows exponentially fast with $n$. It is known, using for instance the technology of relative train tracks (see~\cite{BesHan92}) that, otherwise, $[g]$ has polynomial growth under iterates of $\phi$. Let $\mathrm{Poly}(\phi)$ be the set of conjugacy classes of elements of $F_n$ whose growth under iteration of $\phi$ is polynomial. For a subgroup $H$ of $F_n$, let $\mathrm{Poly}(H)=\bigcap_{\phi \in H} \mathrm{Poly}(\phi)$. The aim of these three papers is to prove the following result: 

\begin{theo}[\cite{Guerch2021Polygrowth}]\label{Theo intro 2}
Let $n \geq 3$ and let $H$ be a subgroup of $\Out(F_n)$. There exists $\phi \in H$ such that $\mathrm{Poly}(\phi)=\mathrm{Poly}(H)$.
\end{theo}

Theorem~\ref{Theo intro 2} is proved using dynamical methods developed mainly in~\cite{Guerch2021NorthSouth}. In the present article, we introduce the topological space associated with the dynamics. Informally, Theorem~\ref{Theo intro 2} shows that the exponential growth of a subgroup $H$ of $\Out(F_n)$ is encaptured by the exponential growth of a single element of $H$. In this paper, we construct a space which is well-adapted for our considerations, the \emph{space of currents relative to a malnormal subgroup system}. These relative currents are positive $F_n$-invariant Radon measures on an appropriate subspace of the double boundary at infinity of $F_n$. Let $\phi \in \Out(F_n)$. When the malnormal subgroup system is appropriately chosen, this space has the property that its points corresponding to conjugacy classes of elements in $F_n-\mathrm{Poly}(\phi)$ are dense in it (see~Theorem~\ref{Theo introduction}).

The space of currents that we construct in this paper builds on objects introduced for similar purposes. For instance, the study of the mapping class group $\Mod(S)$ of a connected, compact, oriented surface $S$ has benefited from the study of the action of $\Mod(S)$ on the \emph{space of geodesic currents} $\Curr(S)$, introduced by Ruelle and Sullivan in~\cite{RuelleSullivan1975} (see also the work of Bonahon~\cite{Bonahon86}). It is defined as the space of $\pi_1(S)$-invariant and flip invariant nonnegative Radon measures on the double boundary $\partial^2 \widetilde{S}$ of a universal cover $\widetilde{S}$ of $S$, equipped with the weak-star topology. Considering the space of projective geodesic currents $\PCurr(S)$, one can show that $\PCurr(S)$ can be viewed as a completion of the currents associated with weighted nontrivial homotopy classes of closed curves on $S$. The space $\PCurr(S)$ is well-adapted to the study of $\Mod(S)$. For instance, it can be used for counting closed geodesics whose length is bounded by a given constant when the surface $S$ is equipped with a hyperbolic metric (see~\cite{ErlanUya2020} for a survey). Concerning dynamical properties, a result of Thurston (\cite{Thurston88}, see also~\cite{Uyanik2015}) implies that pseudo-Anosov diffeomorphisms act with North-South dynamics on the space $\PCurr(S)$: every pseudo-Anosov element $f \in \Mod(S)$ has exactly two fixed points in $\PCurr(S)$ and any other nonfixed point in $\PCurr(S)$ converges to one of the fixed points under positive or negative iterates of $f$. Moreover, this convergence can be made uniform on compact subsets of $\PCurr(S)$ which do not contain the fixed points. 

In the specific context of free groups, building on~\cite{Bonahon91} for general hyperbolic groups, the space of currents $\Curr(F_n)$ was first studied by Martin~\cite{Martin95}. It is defined as the space of $F_n$-invariant, flip invariant nonnegative Radon measure on the double boundary $\partial^2F_n$ of $F_n$ equipped with the weak-star topology. Martin showed that the set of currents associated with conjugacy classes of nontrivial elements of $F_n$ is dense in the space $\PCurr(F_n)$ of projective currents. Currents for free groups have also been studied in~\cite{Kapovich2006,KapLus09,CouHilLus08III}. Similarly to pseudo-Anosov elements of $\Mod(S)$ on $\PCurr(S)$, fully irreducible automorphisms of $F_n$ and atoroidal automorphisms of $F_n$ act with North-South type  dynamics on $\PCurr(F_n)$ (see~\cite{Uyanik2015,Uyanik2019}).

Currents on free groups have also been studied in a relative context, more precisely, in the context of \emph{free factor systems}. A free factor system $\mathcal{F}$ is a finite set of conjugacy classes $\mathcal{F}=\{[A_1],\ldots,[A_k]\}$ of nontrivial subgroups $A_1,\ldots,A_k$ of $F_n$ such that there exists a subgroup $B$ of $F_n$ with $F_n=A_1 \ast \ldots A_k \ast B$. Gupta~\cite{gupta2017relative} (see also~Guirardel-Horbez~\cite{Guirardelhorbez19laminations}) introduced the space $\Curr(F_n,\mathcal{F})$ of currents relative to the free factor system $\mathcal{F}$. Relative currents are then $F_n$-invariant, flip invariant nonnegative Radon measures on a subspace of the double boundary of $F_n$ which does not intersect the double boundary of any conjugate of $A_i$, equipped with the weak-star topology. Gupta~\cite{gupta2017relative} then showed that the set of currents associated with conjugacy classes of \emph{nonperipheral elements of $F_n$}, that is, elements of $F_n$ that do not belong to any conjugate of some $A_i$, is dense in $\PCurr(F_n,\mathcal{F})$. She then showed that fully irreducible outer automorphisms relative to $\mathcal{F}$ act with a North-South type dynamics on $\PCurr(F_n,\mathcal{F})$.

In order to study the purely exponential growth part of an outer automorphism of $F_n$, we need to consider currents relative to a class of subgroup systems which is larger than the class of free factor systems. Indeed, if $\phi \in \Out(F_n)$, the set of all maximal conjugacy classes of subgroups of $F_n$ consisting of elements with polynomial growth under iterates of $\phi$ is not necessarily a free factor system. However, Levitt~\cite[Proposition~1.4]{Levitt09} proved that this set is a \emph{malnormal subgroup system}. A malnormal subgroup system $\mathcal{A}$ is a finite set of conjugacy classes $\mathcal{A}=\{[A_1],\ldots, [A_k]\}$ of nontrivial subgroups of $F_n$ such that, for every $i \in \{1,\ldots,k\}$, the group $A_i$ is malnormal and, for every subgroups $B_1,B_2$ of $F_n$ such that $[B_1],[B_2] \in \mathcal{A}$, if the intersection $B_1 \cap B_2$ is nontrivial, then $B_1=B_2$. A free factor system is, in particular, a malnormal subgroup system but the converse does not hold (see~Section~\ref{Section malnormal subgroup}).

Let $\mathcal{A}=\{[A_1],\ldots, [A_k]\}$ be a malnormal subgroup system. We define the space $\Curr(F_n,\mathcal{A})$ of \emph{currents relative to $\mathcal{A}$} as the space of $F_n$-invariant, flip invariant nonnegative Radon measures on a natural space $\partial^2(F_n,\mathcal{A})$, the double boundary of $F_n$ relative to $\mathcal{A}$ equipped with the weak-star topology. The space $\partial^2(F_n,\mathcal{A})$ is a subspace of $\partial^2F_n$ which does not intersect the double boundary of any conjugate of $A_i$ (see Section~\ref{Section double boundary} for precise definitions). In this article, we prove the following result. An element of $F_n$ is \emph{non-$\mathcal{A}$-peripheral} if it is not contained in any conjugate of any $A_i$ with $i \in \{1,\ldots,k\}$.

\begin{theo}\label{Theo introduction}
Let $n \geq 3$ and let $\mathcal{A}$ be a malnormal subgroup system. The set of currents associated with conjugacy classes of non-$\mathcal{A}$-peripheral elements of $F_n$ is dense in the space $\PCurr(F_n,\mathcal{A})$ of projective currents relative to $\mathcal{A}$.
\end{theo}

Let $\phi \in \Out(F_n)$. If $\mathcal{A}$ is the set of conjugacy classes of maximal polynomial subgroups of $\phi$, then Theorem~\ref{Theo introduction} shows that the set of projective currents associated with exponentially growing elements of $F_n$ under iterates of $\phi$ is dense in $\PCurr(F_n,\mathcal{A})$. Therefore, the space $\PCurr(F_n,\mathcal{A})$ is a natural topological space for the study of the action of $\phi$ on elements of $F_n$ with exponential growth under iterates of $\phi$. A subsequent paper~\cite{Guerch2021NorthSouth} will then show that $\phi$ acts with North-South type dynamics on $\PCurr(F_n,\mathcal{A})$. This North-South dynamics will be a central argument in the proof of Theorem~\ref{Theo intro 2}.

We now give an outline of the proof of Theorem~\ref{Theo introduction}. The proof follows the one of a similar result in the context of currents relative to free factor systems due to Gupta~\cite{gupta2017relative}. However, in the case of free factor systems, the proof relies on the existence of an adapted free basis of $F_n$ associated with the free factor system, which does not necessarily exist in the case of malnormal subgroup systems. Our new argument in order to overcome this difficulty is the description of a finite set of elements of $F_n$ associated with a malnormal subgroup system and a free basis of $F_n$ which completely determines whether an element of $F_n$ is contained in a conjugate of a subgroup of the malnormal subgroup system or not (see Lemma~\ref{Lem finite set of words determines vertex group system}). 

Let $\mathcal{A}$ be a malnormal subgroup system and let $\mu \in \PCurr(F_n,\mathcal{A})$. We first show that $\mu$ can be extended into a \emph{signed measured current $\widetilde{\mu}$ on $F_n$}, that is  an $F_n$-invariant and flip invariant Radon measure on $\partial^2F_n$. Even though $\widetilde{\mu}$ might have negative values, we show that $\widetilde{\mu}$ can be chosen so that $\widetilde{\mu}$ gives positive value to sufficiently many Borel subsets of $\partial^2 F_n$. One can then use the density of currents associated with conjugacy classes of nontrivial elements of $F_n$ in the space $\Curr(F_n)$ in order to conclude the proof.

To our knowledge, the objects we construct in this paper have not been studied or constructed for larger classes of groups, such as relatively hyperbolic groups and quasi-convex almost malnormal subgroups of hyperbolic groups. Nevertheless, the extension of our definitions to this context seems natural since a result of Bowditch~\cite[Theorem~7.11]{Bowditch2012} shows that the group $F_n$ is always hyperbolic relative to a malnormal subgroup system $\mathcal{A}$. But as we explained in Remark~\ref{Rmq other possibilities for boundary}, the natural double boundary associated with a relative hyperbolic group will have less information than the boundary $\partial^2(F_n,\mathcal{A})$. Therefore, it would require new techniques to develop the notion of currents for relative hyperbolic groups or quasi-convex almost malnormal subgroups of hyperbolic groups.

\medskip

{\small{\bf Acknowledgments. } I warmly thank my advisors, Camille Horbez and Frédéric Paulin, for their precious advices and for carefully reading the different versions of this article.}

\section{Malnormal subgroup systems}\label{Section malnormal subgroup}

\subsection{Malnormal subgroup systems}

Let $n$ be an integer greater than $1$ and let $F_n$ be a free group of rank $n$. In this section, we define, following Handel and Mosher~\cite[Section~I.1.1.2]{HandelMosher20}, malnormal subgroups systems and study some of their properties.

A \emph{subgroup system of $F_n$} is a finite (possibly empty) set $\mathcal{A}$ whose elements are conjugacy classes of nontrivial (that is distinct from $\{1\}$ and $F_n$) finite rank subgroups of $F_n$. Note that a subgroup system $\mathcal{A}$ is completely determined by the set of subgroups $A$ of $F_n$ such that $[A] \in \mathcal{A}$. There exists a preorder on the set of subgroup systems of $F_n$, where $\mathcal{A}_1 \leq \mathcal{A}_2$ if for every subgroup $A_1$ of $F_n$ such that $[A_1] \in \mathcal{A}_1$, there exists a subgroup $A_2$ of $F_n$ such that $[A_2] \in \mathcal{A}_2$ and $A_1$ is a subgroup of $A_2$. The \emph{stabilizer in $\Out(F_n)$ of a subgroup system} $\mathcal{A}$, denoted by $\Out(F_n,\mathcal{A})$, is the set of all elements $\phi \in \Out(F_n)$ such that $\phi(\mathcal{A})=\mathcal{A}$.

Recall that a subgroup $A$ of $F_n$ is \emph{malnormal} if for every element $x \in F_n-A$, we have $xAx^{-1} \cap A=\{e\}$. A subgroup system $\mathcal{A}$ is said to be \emph{malnormal} if every subgroup $A$ of $F_n$ such that $[A] \in \mathcal{A}$ is malnormal and, for any subgroups $A_1,A_2$ of $F_n$ such that $[A_1],[A_2] \in \mathcal{A}$, if $A_1 \cap A_2$ is nontrivial then $A_1=A_2$. There are equivalent formulations of malnormality which we present now (see~\cite[Section~I.1.1.2]{HandelMosher20}). Let $T$ be the Cayley graph of $F_n$ with respect to some given free basis of $F_n$. For every subgroup $A$ of $F_n$, let $T_A$ be the minimal $A$-invariant subtree of $T$. Then a subgroup system $\mathcal{A}$ made of conjugacy classes of malnormal subgroups is malnormal if and only if there exists a finite constant $L >0$ such that for any distinct subgroups $A_1,A_2$ of $F_n$ such that $[A_1],[A_2] \in \mathcal{A}$, the diameter of the intersection $T_{A_1} \cap T_{A_2}$ is at most equal to $L$. Malnormality of a subgroup system $\mathcal{A}$ made of conjugacy classes of malnormal subgroups is also equivalent to the fact that, for any distinct subgroups $A_1$ and $A_2$ of $F_n$ such that $[A_1],[A_2] \in \mathcal{A}$, we have $\partial_{\infty}  T_{A_1} \cap \partial_{\infty} T_{A_2}=\varnothing$. 

\subsection{Properness at infinity}

Let $\partial_{\infty}F_n$ be the Gromov boundary of $F_n$. Let $\mathcal{B}$ be a free basis of $F_n$ and let $T$ be the Cayley graph of $F_n$ with respect to $\mathcal{B}$. For convenience, we suppose that $\mathcal{B}^{-1}=\mathcal{B}$. The boundary of $T$ is naturally homeomorphic to $\partial_{\infty}F_n$. For an element $w \in F_n$, we denote by $\gamma_w$ the path in $T$ starting from $e$ corresponding to the word $w$. We denote by $w^{+\infty}$ the element in $\partial_{\infty}F_n$ corresponding to the quasi-geodesic starting at $e$ obtained by concatenating paths in $T$ labeled by $w$.

Let $A$ be a subgroup of $F_n$ of finite rank. The inclusion $A \subseteq F_n$ induces an $A$-equivariant inclusion $\partial_{\infty} A \hookrightarrow \partial_{\infty} F_n$. Note that the $F_n$-orbit of the image of this map only depends on the conjugacy class of $A$ in $F_n$.

Let $\mathcal{A}$ be a subgroup system of $F_n$. The subgroup system $\mathcal{A}$ is said to be \emph{proper at infinity} if, for every element $g$ of $F_n$, the following assertions are equivalent:

\begin{itemize}
\item there exists a subgroup $A$ of $F_n$ such that $[A] \in \mathcal{A}$ and $g^{+\infty} \in \partial_{\infty} A$;
\item there exists a subgroup $A$ of $F_n$ such that $[A] \in \mathcal{A}$ and $g \in A$.
\end{itemize} 

For the proof of Lemma~\ref{Lem malnormal is proper at infinity} below, we need the following result (see for instance~\cite[Fact~1.2]{HandelMosher20}). This is a particular case of the same result valid for all
quasi-convex subgroups $A_1,A_2$ of any word hyperbolic group, see~\cite{Swenson97}, that has been for instance generalized in \cite[Theorem~1.4]{Tran19}.

\begin{lem}\label{Lem intersection at infinity}
For every finitely generated subgroups $A_1$ and $A_2$ of $F_n$, we have $$\partial_{\infty} (A_1 \cap A_2)=\partial_{\infty} A_1 \cap \partial_{\infty} A_2.$$ 
\end{lem}

A subgroup $A$ of $F_n$ is \emph{root-closed} if for every $g \in F_n$ and every $k \in \NN^*$ such that $g^k \in A$, we have $g \in A$.

\begin{lem}\label{Lem malnormal is proper at infinity}
Let $\mathcal{A}$ be a subgroup system. The following are equivalent:

\medskip

\noindent{$(1)$ } the subgroup system $\mathcal{A}$ is proper at infinity;

\medskip

\noindent{$(2)$ } every subgroup $A$ of $F_n$ such that $[A] \in \mathcal{A}$ is root-closed. 

\medskip

In particular, a malnormal subgroup system is proper at infinity.
\end{lem}

\dem Suppose that $\mathcal{A}$ is proper at infinity and let $A$ be a subgroup of $F_n$ such that $[A] \in \mathcal{A}$. Let $g \in F_n$ and $k \in \NN^*$ be such that $g^k \in A$. Let us prove that $g \in A$. Since $g^k \in A$, we see that $g^{+\infty} \in \partial_{\infty} A$. Since $\mathcal{A}$ is proper at infinity, we have $g \in A$. Hence $A$ is root-closed. Suppose now that every subgroup $A$ of $F_n$ such that $[A] \in \mathcal{A}$ is root-closed. Let $g \in F_n$ and let $A$ be a subgroup of $F_n$ such that $[A] \in \mathcal{A}$ and $g^{+\infty} \in \partial_{\infty} A$. By Lemma~\ref{Lem intersection at infinity} applied to $\left\langle g \right\rangle$ and $A$, there exists $k \in \NN^*$ such that $g^k \in A$. Since $A$ is root-closed, we see that $g \in A$. Hence $\mathcal{A}$ is proper at infinity. This shows the equivalence.

Let $\mathcal{A}$ be a malnormal subgroup system and let $A$ be a subgroup of $F_n$ such that $[A] \in \mathcal{A}$. We prove that $A$ is root-closed. Let $g \in F_n$ and let $k \in \NN^*$ be such that $g^k \in A$. We claim that $g \in A$. Indeed, suppose towards a contradiction that $g \notin A$. Then $g^k=gg^kg^{-1}$ belongs to $A \cap gAg^{-1}$ which is equal to $\{e\}$, a contradiction.
\hfill\qedsymbol

\bigskip

Let $\mathcal{A}$ be a malnormal subgroup system. An element $g \in F_n$ is \emph{$\mathcal{A}$-peripheral} (or simply \emph{peripheral} if there is no ambiguity) if it is trivial or conjugate into one of the subgroups of $\mathcal{A}$, and \emph{non-$\mathcal{A}$-peripheral} otherwise. Note that, since $\mathcal{A} \neq \{[F_n]\}$, there always exists a non-$\mathcal{A}$-peripheral element. Since $\mathcal{A}$ is proper at infinity by Lemma~\ref{Lem malnormal is proper at infinity}, we see that an element $g$ of $F_n$ is $\mathcal{A}$-peripheral if and only if there exists a subgroup $A$ of $F_n$ such that $[A] \in \mathcal{A}$ and $g^{+\infty} \in \partial_{\infty} A$.

Let $\mathcal{A}=\{[A_1],\ldots,[A_r]\}$ be a malnormal subgroup system of $F_n$. For every element \mbox{$i \in \{1,\ldots,r\}$}, let $T_{A_i}$ be the minimal $A_i$-invariant subtree of $T$. Suppose that for every $i \in \{1,\ldots,r\}$, the representative $A_i$ of $[A_i]$ is chosen so that the tree $T_{A_i}$ contains the base point $e$ of $T$.

By malnormality of $\mathcal{A}$, there exists $L \in \NN^*$ such that for any distinct subgroups $A,B$ of $F_n$ such that $[A],[B] \in \mathcal{A}$, the diameter of the intersection $T_{A} \cap T_{B}$ is at most $L$. Let $i \in \{1,\ldots,r\}$. Let $\Gamma_i$ be the set of subgroups $B$ of $F_n$ such that there exists $g_B \in F_n$ such that $B=g_BA_ig_B^{-1}$ and the tree $T_B$ contains the base point $e$ of $T$. Note that, by malnormality of $\mathcal{A}$, for every $i \in \{1,\ldots,r\}$, the set $\Gamma_i$ is finite. Let $C_i$ be the set of elements $w$ of $F_n$ such that the length of $\gamma_w$ is equal to $L+2$ and, for every $B \in \Gamma_i$, the path $\gamma_w$ is not contained in $T_{B}$. Let $\mathscr{C}=\bigcap_{i=1}^r C_i$. Since we are looking at geodesic paths of length equal to $L+2$, the set $\mathscr{C}$ is finite. If $\gamma$ is a path in $T$, the \emph{element of $F_n$ corresponding to $\gamma$} is the element $h \in F_n$ such that the path $\gamma$ is labeled by $h$. 

\begin{lem}\label{Lem finite set of words determines vertex group system}
Let $\mathcal{B}$, $T$, $\mathcal{A}=\{[A_1],\ldots,[A_r]\}$, $L \in \NN^*$, $\Gamma_1,\ldots,\Gamma_r$, $\mathscr{C}$ be as above. The finite set $\mathscr{C}=\mathscr{C}(A_1,\ldots, A_r)$ is nonempty. Moreover, it satisfies the following:

\medskip

\noindent{$(1)$ } every element $g \in F_n$ such that the length of $\gamma_g$ is at least equal to $L+2$ and such that $\gamma_g$ is not contained in a tree $T_B$ with $B \in \bigcup_{i=1}^r \Gamma_i$ contains an element of $\mathscr{C}$ as a subword. In particular, every non-$\mathcal{A}$-peripheral cyclically reduced element $g \in F_n$ has a power which contains an element of $\mathscr{C}$ as a subword;

\medskip

\noindent{$(2)$ } for every non-$\mathcal{A}$-peripheral cyclically reduced element $g \in F_n$, if $c_g$ is the geodesic ray in $T$ starting from $e$ obtained by concatenating edge paths labeled by $g$, there exists an edge path in $c_g$ labeled by a word in $\mathscr{C}$ at distance at most $L+2$ from $\bigcup_{i=1}^r \bigcup_{B \in \Gamma_i} T_B$;

\medskip

\noindent{$(3)$ } if an element $w \in F_n$ contains an element of $\mathscr{C}$ as a subword, then for every $i \in \{1,\ldots,r\}$, the element $w$ is not contained in $A_i$.

\end{lem}

\dem We first prove that $(1)$ and $(2)$ hold and that $\mathscr{C}$ is nonempty. Let $g$ be as in the first claim of Assertion~$(1)$. First note that, by the choice of $L$, for every $i,j \in \{1,\ldots,r\}$ and every distinct $A \in \Gamma_i$ and $B \in \Gamma_j$, the intersection $T_A \cap T_B$ is contained in the closed ball of radius $L$ centered at $e$. We consider the geodesic path $c_g \colon [0,1] \to T$ such that $c(0)=e$ and such that $c_g(1)$ is the terminal endpoint of $\gamma_g$. Let $i \in \{1,\ldots,r\}$ and let $$t_0=\max\left\{t \in [0,1]\;|\; c_g(t) \in \bigcup_{i=1}^r\bigcup_{A \in \Gamma_i} T_{A}\right\}.$$ The point $c_g(t_0)$ is a vertex and is distinct from $c_g(1)$ by assumption. We denote by $c_{\mathcal{A}}$ the geodesic segment $c_g \cap \bigcup_{i=1}^r\bigcup_{A \in \Gamma_i} T_{A}$. 

Suppose first that the length of $c_{\mathcal{A}}$ is at most equal to $L+1$. Let $c_0$ be the geodesic segment contained in $c_g$ which originates at $c_g(t_0)$ and such that the length of $c_{\mathcal{A}}c_0$ is equal to $L+2$. Such a path $c_{\mathcal{A}}c_0$ exists since the length of $\gamma_g$ is at least equal to $L+2$. Then the element $h$ of $F_n$ corresponding to $c_{\mathcal{A}}c_0$ is in $\mathscr{C}$ and is a subword of $g$. This concludes the proof in this case.

Suppose now that the length of $c_{\mathcal{A}}$ is greater than $L+1$. Let $c_{\mathcal{A}}(t_0-L-1)$ be the vertex in $c_{\mathcal{A}}$ at distance $L+1$ from $c_g(t_0)$, and let $g_0$ be the corresponding element of $F_n$. Let $s_0$ be the geodesic path between $c_{\mathcal{A}}(t_0-L-1)$ and $c_g(t_0)$. Since the geodesic path $s_0$ has length equal to $L+1$, there exists a unique $i \in \{1,\ldots,r\}$ and a unique $A \in \Gamma_i$ such that $s_0$ is contained in $T_A$. Let $e_0$ be the edge in $c_g$ which originates at $c_g(t_0)$. Let $h \in F_n$ be the element corresponding to the edge path $s_1$ between $c_g(t_0-L-1)$ and the terminal point of $e_0$. We claim that $h \in \mathscr{C}$. Indeed, suppose towards a contradiction that $h \notin \mathscr{C}$. Then there exists $j \in \{1,\ldots,r\}$ and $B \in \Gamma_j$ such that the edge path $\gamma_{h}$ is contained in $T_{B}$. Since $\gamma_h$ has length equal to $L+2$, the integer $j$ and the subgroup $B$ are unique. Remark that $g_0^{-1}$ sends the geodesic path $s_0$ to the initial segment of length $L+1$ of $\gamma_h$. Since $g_0^{-1}s_0$ has length equal to $L+1$, the subgroup $B$ is the unique element of $\bigcup_{\ell=1}^r \Gamma_{\ell}$ such that the tree $T_B$ contains $g_0^{-1}s_0$. But $s_0$ is contained in $T_{A}$ and the tree $T_{A}$ is sent by $g_0^{-1}$ to the tree $T_{g_0^{-1}Ag_0}$. Therefore, we see that $B=g_0^{-1}Ag_0$. But $g_0^{-1}$ induces an isometry between $T_A$ and $T_{g_0^{-1}Ag_0}$. Therefore, since $s_1$ is not contained in $T_A$, we see that $\gamma_h=g_0^{-1}s_1$ is not contained in $T_{g_0^{-1}Ag_0}$. This leads to a contradiction. Hence $h \in \mathscr{C}$ and $h$ is a subword of $g$. This proves the first claim of Assertion $(1)$. We now prove the second claim of Assertion~$(1)$. Let $g$ be a non-$\mathcal{A}$-peripheral cyclically reduced element of $F_n$. Let $c_g' \colon \RR_+ \to T$ be the geodesic ray in $T$ starting from $e$ obtained by concatenating edge paths labeled by $g$. Recall that, for every $i \in \{1,\ldots,r\}$, the set $\Gamma_i$ is finite. Therefore, since $g$ is nonperipheral and since $\mathcal{A}$ is proper at infinity by Lemma~\ref{Lem malnormal is proper at infinity}, the intersection of $c_g'$ with $\bigcup_{i=1}^r\bigcup_{A \in \Gamma_i} T_{A}$ is compact. Hence there exists a power of $g$ which satisfies the first claim of Assertion~$(1)$. This proves~$(1)$. Moreover, the terminal endpoint of the path in $c_g$ labeled by $h$ which we have constructed is either at distance $L+2$ from $e$ or is at distance at most $1$ from $\bigcup_{i=1}^r \bigcup_{B \in \Gamma_i} T_B$. This proves $(2)$. This also proves that $\mathscr{C}$ is nonempty as there exists a non-$\mathcal{A}$-peripheral element.

We now prove $(3)$. Suppose towards a contradiction that there exist $i \in \{1,\ldots,r\}$ and $a \in A_i$ such that $a$ contains a word of $\mathscr{C}$ as a subword. Thus there exist $x \in \mathscr{C}$, $b,c \in F_n$ such that $a=bxc$ and the word $bxc$ is reduced. Then since $e$ is contained in $T_{A_i}$, the path $\gamma_a$ is contained in $T_{A_i}$. But the element $b^{-1}$ sends the tree $T_{A_i}$ to the tree $T_{b^{-1}A_ib}$. Moreover, since $T_{A_i}$ contains the vertex labeled by $b$, the tree $T_{b^{-1}A_ib}$ contains the base point $e$ of $T$. But then $T_{b^{-1}A_ib}$ contains the geodesic segment $\gamma_x$. This contradicts the fact that $x \in \mathscr{C} \subseteq C_i$. This concludes the proof.
\hfill\qedsymbol

\subsection{Examples of malnormal subgroup systems}

Let $n$ be an integer greater than $1$ and let $F_n$ be a free group of rank $n$. In this section, we give some examples of malnormal subgroup systems. The first one that we describe, following Handel and Mosher \cite{HandelMosher20}, is an \emph{$\RR$-vertex group system}. Let $T$ be an $\RR$-tree equipped with a minimal, isometric action of $F_n$ for which no point or end of $T$ is fixed by the whole group and with trivial arc stabilizers. A proper, nontrivial subgroup $A$ of $F_n$ is an \emph{$\RR$-vertex group of $T$} if there exists a point $x \in T$ such that $A=\Stab(x)$. Note that every free factor of $F_n$ is an $\RR$-vertex group of some simplicial tree. Every $\RR$-vertex group has rank at most equal to $n$ (see~\cite{GabLev95}). 

The \emph{$\RR$-vertex group system of $T$}, denoted by $\mathcal{A}_T$, is the set consisting of all conjugacy classes of nontrivial point stabilizers in $T$. The set $\mathcal{A}_T$ is finite and its cardinality is bounded from above by a finite constant depending only on $n$ (see~\cite{GabLev95}). Therefore the set $\mathcal{A}_T$ is a subgroup system. Note that every free factor system of $F_n$ is an $\RR$-vertex group system of some simplicial tree. However, there exist $\RR$-vertex group systems which are not free factor systems. For example, let $S$ be a compact connected oriented hyperbolic surface with one totally geodesic boundary component such that $\pi_1(S)$ is isomorphic to $F_n$. Let $T$ be the $\RR$-tree dual to the lift $\widetilde{\Lambda}$ to $\mathbb{H}_2$ of a measured geodesic lamination $\Lambda$ without compact leaves on $S$. An identification of $\pi_1(S)$ with $F_n$ induces an action of $F_n$ on $T$ which has trivial arc stabilizers. Moreover, the fundamental group of the connected component containing the boundary curve of $S$ is the stabilizer of a point in $T$.
Since the fundamental group of this connected component is not a free factor of $F_n$, this shows that $\mathcal{A}_T$ is not a free factor system. More generally, Handel and Mosher \cite[Proposition~3.3]{HandelMosher20} give general constructions of $\RR$-vertex group systems which are not free factor systems. 

\begin{lem}\cite[Lemma~3.1]{HandelMosher20}\label{Lem properties of vertex group system}
The subgroup system $\mathcal{A}_T$ is a malnormal subgroup system.
\end{lem}

Another example of malnormal subgroup systems is the following. An outer automorphism $\phi \in \Out(F_n)$ is \emph{exponentially growing} if there exists $g \in F_n$ such that the length of the conjugacy class $[g]$ of $g$ in $F_n$ with respect to some basis of $F_n$ grows exponentially fast under iteration of $\phi$. If $\phi \in \Out(F_n)$ is not exponentially growing, then the length of the conjugacy class of every element of $F_n$ is polynomially growing under iteration of $\phi$  and $\phi$ is said to be \emph{polynomially growing}. One similarly says that an automorphism $\alpha \in \Aut(F_n)$ is exponentially growing or polynomially growing. Let $\phi \in \Out(F_n)$ be exponentially growing. A subgroup $P$ of $F_n$ is a \emph{polynomial subgroup} of $\phi$ if there exist $k \in \NN^*$ and a representative $\alpha$ of $\phi^k$ such that $\alpha(P)=P$ and $\alpha|_P$ is polynomially growing. By~\cite[Proposition~1.4]{Levitt09}, there exist finitely many conjugacy classes $[H_1],\ldots,[H_k]$ of maximal polynomial subgroups of $\phi$ and the set $\mathcal{H}=\{[H_1],\ldots,[H_k]\}$ is a malnormal subgroup system.

\subsection{Double boundary of $F_n$ relative to a malnormal subgroup system}\label{Section double boundary}

In this section, we construct a boundary of $F_n$ relative to a malnormal subgroup system. We follow a similar construction made by Gupta in~\cite[Section~3.1]{gupta2017relative} in the case of the boundary relative to a free factor system.

The \emph{double boundary of $F_n$} is the quotient topological space $$\partial^2F_n=\left(\partial_{\infty} F_n \times \partial_{\infty} F_n \setminus \Delta \right)/\sim,$$ where $\sim$ is the equivalence relation generated by the flip relation $(x,y)\sim(y,x)$ and $\Delta$ is the diagonal, endowed with the diagonal action of $F_n$. We denote by $\{x,y\}$ the equivalence class of $(x,y)$.

Let $\mathcal{A}=\{[A_1],\ldots,[A_r]\}$ be a malnormal subgroup system of $F_n$. Let $\mathcal{B}$, $T$, $L \in \NN^*$, $\Gamma_1,\ldots,\Gamma_r$, $\mathscr{C}$ be as above Lemma~\ref{Lem finite set of words determines vertex group system}. The boundary of $T$ is naturally homeomorphic to $\partial_{\infty}F_n$ and the set $\partial^2F_n$ is then identified with the set of unoriented bi-infinite geodesics in $T$. Let $\gamma$ be a finite geodesic path in $T$. The path $\gamma$ determines a subset in $\partial^2F_n$ called the \emph{cylinder set of $\gamma$}, denoted by $C(\gamma)$, which consists in all unoriented bi-infinite geodesics in $T$ that contain $\gamma$. Such cylinder sets form a basis for a topology on $\partial^2 F_n$, and in this topology, the cylinder sets are both open and compact, hence closed since $\partial^2F_n$ is Hausdorff. The action of $F_n$ on $\partial^2F_n$ has a dense orbit.

Let $A$ be a nontrivial subgroup of $F_n$ of finite rank. The induced $A$-equivariant inclusion $\partial_{\infty} A \hookrightarrow \partial_{\infty} F_n$ induces an inclusion $\partial^2 A \hookrightarrow \partial^2 F_n$. Let $$\partial^2\mathcal{A}= \bigcup_{i=1}^r \bigcup_{g \in F_n} \partial^2 gA_ig^{-1}.$$ Let $\partial^2(F_n,\mathcal{A})=\partial^2F_n -\partial^2\mathcal{A}$ be the \emph{double boundary of $F_n$ relative to $\mathcal{A}$}. This subset is invariant under the action of $F_n$ on $\partial^2F_n$ and inherits the subspace topology of $\partial^2F_n$, denoted by $\tau$.

\begin{lem}\label{Lem double boundary open}
Let $\mathrm{Cyl}(\mathscr{C})$ be the set of cylinder sets of the form $C(\gamma)$, where the element of $F_n$ determined by the geodesic edge path $\gamma$ contains an element of $\mathscr{C}$ as a subword. We have $$\partial^2(F_n,\mathcal{A})=\bigcup_{C(\gamma) \in \mathrm{Cyl}(\mathscr{C})}C(\gamma).$$ In particular, the space $\partial^2(F_n,\mathcal{A})$ is an open subset of $\partial^2 F_n$.
\end{lem}

\dem Let $y \in \partial^2(F_n,\mathcal{A})$. Let $c$ be an oriented geodesic line $c$ in $T$ which belongs to the equivalence class $y$. Let $v$ be a vertex of $T$ contained in $c$ and let $g_0$ be the corresponding element of $F_n$. 

Suppose first that the intersection $c \cap g_0 \left(\bigcup_{i=1}^r \bigcup_{B \in \Gamma_i} T_B\right)$ is either compact or a half-line. In particular, the intersection $c \cap g_0 \left(\bigcup_{i=1}^r \bigcup_{B \in \Gamma_i} T_B\right)$ has a terminal point $v'$. Let $x$ be the vertex in $c$ at distance $L+2$ from $v'$. Let $g \in F_n$ be the element corresponding to the geodesic edge path between $v$ and $x$. Note that the edge path $\gamma_g$ is not contained in $\bigcup_{i=1}^r \bigcup_{B \in \Gamma_i} T_B$ since, for every nontrivial subgroup $A$ of $F_n$ of finite rank, the element $g_0$ sends $T_A$ to $T_{g_0Ag_0^{-1}}$. By Lemma~\ref{Lem finite set of words determines vertex group system}~$(2)$, the word $g$ contains a word of $\mathscr{C}$ as a subword. Then $y \in g_0C(\gamma_g)$, and $g_0C(\gamma_g) \in \mathrm{Cyl}(\mathscr{C})$.

Suppose now that the intersection $c \cap g_0 \left(\bigcup_{i=1}^r \bigcup_{B \in \Gamma_i} T_B\right)$ is not compact. Since $y \in \partial^2(F_n,\mathcal{A})$, the path $c$ cannot be contained in a single tree $g_0T_B$ with $B \in \bigcup_{i=1}^r\Gamma_i$. By the definition of $L$, there exist exactly two subgroups $A,B \in \bigcup_{i=1}^r\Gamma_i$ such that $c$ is contained in $g_0T_A \cup g_0T_B$. By the definition of the constant $L$, the intersection $g_0T_A \cap g_0T_B$ has diameter at most equal to $L$. Let $c_0$ be the subpath of $c$ of length $2L+2$ whose middle point is $v$ and whose starting point is in $g_0T_A$ and let $g$ be the element of $F_n$ corresponding to $c_0$. Let $v'$ be the initial vertex of $c_0$ and let $g'$ be the element of $F_n$ associated with $v'$. Note that the intersection of $c_0$ with $g_0T_A$ and $g_0T_B$ has length at least equal to $L+1$. Up to considering a larger path $c_0$, we may suppose that $g$ is cyclically reduced. We claim that $g$ contains an element of $\mathscr{C}$ as a subword. Indeed, suppose towards a contradiction that $g$ does not contain an element of $\mathscr{C}$ as a subword. By Lemma~\ref{Lem finite set of words determines vertex group system}~$(1)$, there exist $i \in \{1,\ldots,r\}$ and $H \in \Gamma_i$ such that $\gamma_g \subseteq T_H$. But then $g'\gamma_g=c_0$ and is contained in $g'T_H$. Thus the diameter of the intersection $g'T_H$ with $g_0T_A$ and $g_0T_B$ is at least equal to $L+1$. By definition of $L$, this means that $g'T_H=g_0T_A=g_0T_B$. This means that $A=B$, a contradiction. Hence $g$ contains an element of $\mathscr{C}$ as a subword. Thus we have $y \in g_0C(\gamma_g)$, with $g_0C(\gamma_g) \in \mathrm{Cyl}(\mathscr{C})$. Therefore, we see that $$\partial^2(F_n,\mathcal{A})\subseteq \bigcup_{C(\gamma) \in \mathrm{Cyl}(\mathscr{C})}C(\gamma).$$ 

Conversely, let $\gamma$ be a geodesic path in $T$ such that $C(\gamma) \in \mathrm{Cyl}(\mathscr{C})$. Suppose towards a contradiction that there exists $y \in \partial^2 \mathcal{A}$ such that $y \in C(\gamma)$. Thus, there exist elements $i \in \{1,\ldots,r\}$, $g \in F_n$ and $a \in gA_ig^{-1}$ such that $\{a^{+\infty},a^{-\infty}\} \in C(\gamma)$. Therefore, we see that $\gamma$ is a subpath of $T_{gA_ig^{-1}}$. Decompose $\gamma$ as $\gamma=\tau_1\delta\tau_2$ where $\delta$ is labeled by a word $w$ in $\mathscr{C}$. Let $v$ be the origin of $\delta$ and let $h$ be the element of $F_n$ corresponding to $v$. Then $h^{-1}T_{gA_ig^{-1}}=T_{h^{-1}gA_ig^{-1}h} \in \Gamma_i$ and contains $\gamma_w$ with $w \in \mathscr{C}$, a contradiction.
\hfill\qedsymbol

\bigskip

Note that Lemma~\ref{Lem double boundary open} implies that we can define a topology on $\partial^2(F_n,\mathcal{A})$, denoted by $\tau'$, where cylinder sets in $\mathrm{Cyl}(\mathscr{C})$ generate the topology. Lemma~\ref{Lem double boundary open} also implies that the two topologies $\tau$ and $\tau'$ are equal. 

Since $\partial^2 F_n$ is locally compact and since $\partial^2(F_n,\mathcal{A})$ is an open subset of $\partial^2F_n$ by Lemma~\ref{Lem double boundary open}, we have the following result.

\begin{lem}\label{Lem boundary loc compact}
The space $\partial^2(F_n,\mathcal{A})$ is locally compact.
\hfill\qedsymbol
\end{lem}

\begin{lem}\label{Lem boundary cocompact action}
The action of $F_n$ on $\partial^2(F_n,\mathcal{A})$ has a dense orbit.
\end{lem}

\dem Recall that there exists $g \in F_n$ such that the action of $g$ on $\partial^2F_n$ has a dense orbit. Since $\partial^2(F_n,\mathcal{A})$ is an open subset of $\partial^2F_n$, the element $g$ also acts on $\partial^2(F_n,\mathcal{A})$ with a dense orbit.
\hfill\qedsymbol

\begin{rmq}\label{Rmq other possibilities for boundary}
We now compare our definition with other natural constructions of double boundaries. The first one is to see the double boundary of $F_n$ relative to a malnormal subgroup system as the double boundary of a Gromov hyperbolic space. Indeed, if $\mathcal{A}=\{[A_1],\ldots,[A_r]\}$ is a malnormal subgroup system, by a result of Bowditch (see~\cite[Theorem~7.11]{Bowditch2012}), the group $F_n$ is hyperbolic relative to $\mathcal{A}$. In particular, there is a natural (that is well-defined up to quasi-isometry) proper geodesic Gromov hyperbolic space $X$ on which $F_n$ acts by isometries and such that the subgroups of $F_n$ whose conjugacy classes are in $\mathcal{A}$ are precisely the maximal parabolic subgroups of the action of $F_n$ on the Gromov-boundary of $X$ (see~\cite{Bowditch2012} for a precise description of $X$). Thus a natural construction for another type of double boundary of $F_n$ relative to $\mathcal{A}$ is to define it as the double boundary of $X$. This definition seems to extend to the more general case of relatively hyperbolic groups. However, the relative double boundary $\partial^2(F_n,\mathcal{A})$ has the advantage of being an open subset of $\partial^2 F_n$, so that one can use the cylinder sets of $\partial^2 F_n$ as a basis for the topology of $\partial^2(F_n,\mathcal{A})$. Moreover, the natural application from $\partial F_n$ to $\partial X$ sends the boundary of a parabolic subgroup to a point. Therefore, the relative double boundary $\partial^2(F_n,\mathcal{A})$ seems to contain more information about the geodesic lines whose endpoints are in the Gromov boundary of distinct parabolic subgroups. 

Another candidate for the double boundary of the pair $(F_n,\mathcal{A})$ is the following. Let $\widehat{T}$ be the graph obtained from $T$ by adding one vertex $v(gA)$ for every coset $gA$ with $A$ a subgroup of $F_n$ such that $[A] \in \mathcal{A}$ and by adding an edge from $v(gA)$ to every vertex of $T$ labeled by an element in $gA$. The graph $\widehat{T}$ is Gromov hyperbolic (see for instance~\cite[Proposition~2.6]{KapovichRafi14} or~\cite{Bowditch2012}) and the Gromov boundary of $\widehat{T}$ is homeomorphic to the space $\partial_{\infty} F_n-\bigcup_{i=1}^r\bigcup_{g \in F_n} \partial_{\infty} gA_i$ (see for instance~\cite[Theorem~1.6]{AbbottManning2021} or~\cite{DowdallTaylor2018,Bowditch2012}). However, the double boundary $\partial^2 \widehat{T}$ does not contain any geodesic line whose endpoints are in distinct parabolic subgroups, which makes it a proper subspace of $\partial^2(F_n,\mathcal{A})$ which does not seem to be a union of cylinder sets.
\end{rmq}

\section{Currents relative to a malnormal subgroup system}

In this section, we define currents of $F_n$ relative to a malnormal subgroup system. We follow the construction of Gupta~\cite[Section~3.2]{gupta2017relative} of currents relative to a free factor system.

Let $\mathcal{A}=\{[A_1],\ldots,[A_r]\}$ be a malnormal subgroup system of $F_n$. Let $\mathcal{B}$, $T$, $L \in \NN^*$, $\Gamma_1,\ldots,\Gamma_r$, $\mathscr{C}$ be as above Lemma~\ref{Lem finite set of words determines vertex group system}. 

A \emph{relative current on $(F_n,\mathcal{A})$} is an $F_n$-invariant nonnegative Radon measure $\mu$ on the locally compact space (by Lemma~\ref{Lem boundary loc compact}) $\partial^2(F_n,\mathcal{A})$ (that is $\mu$ gives finite measure to compact subsets of $\partial^2(F_n,\mathcal{A})$, is inner and outer regular). The set $\Curr(F_n,\mathcal{A})$ of all 
relative currents on $\partial^2(F_n,\mathcal{A})$ is equipped with the weak-star topology: a sequence $(\mu_n)_{n \in \NN}$ in $\Curr(F_n,\mathcal{A})^{\NN}$ converges to a current $\mu \in \Curr(F_n,\mathcal{A})$ if and only if for every disjoint clopen subsets $S,S' \subseteq \partial^2(F_n,\mathcal{A})$, the sequence $(\mu_n(S \times S'))_{n \in \NN}$ converges to $\mu(S \times S')$. The space $\Curr(F_n,\mathcal{A})$ is naturally identified with the space of non-negative, $F_n$-invariant, continuous linear functionals on the space $C_c(\partial^2(F_n,\mathcal{A}))$ (equipped with the uniform norm) of continuous compactly supported functions of $\partial^2(F_n,\mathcal{A})$ (see~\cite[Theorem~7.5.5]{Cohn80}). Therefore, the space $\Curr(F_n,\mathcal{A})$ is homeomorphic to a subspace of $C_c(\partial^2(F_n,\mathcal{A}))^*$ equipped with the weak-star topology. Equipped with the uniform structure induced by the weak-star topology on $C_c(\partial^2(F_n,\mathcal{A}))^*$, we see that the space $\Curr(F_n,\mathcal{A})$ is metrisable and complete (see~\cite[Chap.~3, Section~1, Proposition~14]{Bourbaki65}).

The group $\Out(F_n,\mathcal{A})$ acts on $\Curr(F_n,\mathcal{A})$ as follows. Let $\phi \in \Out(F_n,\mathcal{A})$, let $\Phi$ be a representative of $\phi$, let $\mu \in \Curr(F_n,\mathcal{A})$ and let $C$ be a Borel subset of $\partial^2(F_n,\mathcal{A})$. Then, since $\phi$ preserves $\mathcal{A}$, we see that $\Phi^{-1}(C)$ is a Borel subset of $\partial^2(F_n,\mathcal{A})$. Then we set $$\phi(\mu)(C)=\mu(\Phi^{-1}(C)),$$ which is independent of the choice of the representative $\Phi$ since $\mu$ is $F_n$-invariant and the extension to the boundary of the action by conjugation and by left translation of $F_n$ on itself coincide.

We now describe some coordinates for $\Curr(F_n,\mathcal{A})$. Recall that $\mathrm{Cyl}({\mathscr{C}})$ is the set of cylinder sets of the form $C(\gamma)$, where the element of $F_n$ determined by the geodesic path $\gamma$ contains an element of $\mathscr{C}$ as a subword. Recall that $$\partial^2(F_n,\mathcal{A})=\bigcup_{C(\gamma) \in \mathrm{Cyl}(\mathscr{C})}C(\gamma).$$ Let $\eta \in \Curr(F_n,\mathcal{A})$. Let $w \in F_n$ be such that $C(\gamma_w) \in \mathrm{Cyl}({\mathscr{C}})$ and let $w=w_1\ldots w_k$ be the reduced word associated with $w$ written in the basis $\mathcal{B}$. Then $C(\gamma_w)=\coprod C(\gamma_{wb})$, where the union is taken over all elements $b$ of $\mathcal{B} = \mathcal{B}^{-1}$ except $b=w_k^{-1}$. The $\sigma$-additivity of a relative current $\eta$ implies that:
$$\eta(C(\gamma_w))=\sum_{b \neq w_k^{-1}}\eta(C(\gamma_{wb})).$$

Finally, we note that, for every element $w \in F_n$ such that $C(\gamma_w) \in \mathrm{Cyl}({\mathscr{C}})$, we have $\eta(C(\gamma_w))=\eta(C(\gamma_{w^{-1}}))$. Indeed, this follows from the fact that $C(\gamma_w)=wC(\gamma_{w^{-1}})$ and from the $F_n$-invariance of $\eta$. 

\begin{lem}\label{Lem compact subset of boundary disjoint union}
Let $n \geq 3$ and let $C$ be a compact open subset of $\partial^2 (F_n,\mathcal{A})$. There exist finite geodesic edge paths $\gamma_1,\ldots,\gamma_k$ such that:

\medskip

\noindent{$(1)$ } For every $i \in \{1,\ldots,k\}$, we have $C(\gamma_i) \in \mathrm{Cyl}({\mathscr{C}})$;

\medskip

\noindent{$(2)$ } for every distinct $i,j \in \{1,\ldots,k\}$ we have $C(\gamma_i) \cap C(\gamma_j)=\varnothing$;

\medskip

\noindent{$(3)$ } we have $C=\bigcup_{i=1}^k C(\gamma_i)$.
\end{lem}

\dem Since $C$ is a compact open subset of $\partial^2F_n$, using the topology $\tau'$, the set $C$ can be written as a union of cylinder sets $C(\gamma_{1}),\ldots,C(\gamma_{\ell})$, where, for every $i \in \{1,\ldots,\ell\}$, we have $C(\gamma_i) \in \mathrm{Cyl}({\mathscr{C}})$. We may suppose that for every distinct $i,j \in \{1,\ldots,\ell\}$, we have $C(\gamma_i) \nsubseteq C(\gamma_j)$. In particular, there does not exist $i,j \in \{1,\ldots,\ell\}$ such that $\gamma_i \subseteq \gamma_j$. Let $m$ be the number of pairs of distinct elements $i,j \in \{1,\ldots,\ell\}$ such that $C(\gamma_i) \cap C(\gamma_j) \neq \varnothing$. We prove Lemma~\ref{Lem compact subset of boundary disjoint union} by induction on $m$. If for every distinct $i,j \in \{1,\ldots,\ell\}$, we have $C(\gamma_i) \cap C(\gamma_j)=\varnothing$, then the set $\{\gamma_1,\ldots,\gamma_{\ell}\}$ satisfies the conclusion of the lemma. Suppose that there exists $m$ pairs of distinct elements $i,j \in \{1,\ldots,\ell\}$ such that $C(\gamma_i) \cap C(\gamma_j) \neq \varnothing$, with $m \geq 1$.

\medskip

\noindent{\bf Claim. } Let $i,j$ be as above. There exists finite geodesic paths $\gamma_1^{(i)},\ldots,\gamma_{k_i}^{(i)},\gamma_1^{(j)},\ldots,\gamma_{k_j}^{(j)}$ in $T$ which satisfy the following:

\medskip

\noindent{$(a)$ } for every $s \in \{1,\ldots,k_i\}$ and every $t \in \{1,\ldots,k_j\}$, we have $\gamma_i \subseteq \gamma_s^{(i)}$ and $\gamma_j \subseteq \gamma_t^{(j)}$;

\medskip

\noindent{$(b)$ } for every $p \in \{i,j\}$, for every distinct $s,t \in \{1,\ldots, k_p\}$, we have $C(\gamma_s^{(p)}) \cap C(\gamma_t^{(p)})=\varnothing$;

\medskip

\noindent{$(c)$ } for every $s \in \{1,\ldots,k_i\}$ and every $t \in \{1,\ldots,k_j\}$, either $C(\gamma_s^{(i)}) = C(\gamma_t^{(j)})$ or $C(\gamma_s^{(i)}) \cap C(\gamma_t^{(j)})=\varnothing$;

\medskip

\noindent{$(d)$ } for every $p \in \{i,j\}$, we have $$C(\gamma_p)=\bigcup_{s=1}^{k_p} C(\gamma_s^{(p)}).$$

\begin{figure}[ht]
\centering
\begin{tikzpicture}
\draw (-4,0)--(4,0);
\draw (-4,0) node {$\bullet$};
\draw (4,0) node {$\bullet$};
\draw (-3,0) node {$\bullet$};
\draw (3,0) node {$\bullet$};
\draw (-1,0) node {$\bullet$};
\draw (1,0) node {$\bullet$};
\draw (-4,0) node[above] {$v_i'$};
\draw (4,0) node[above] {$v_j'$};
\draw (-1,0) node[above] {$v_i$};
\draw (1,0.01) node[above] {$v_j$};
\draw (-3.5,0) node[below] {$e_i'$};
\draw (3.5,0) node[below] {$e_j'$};
\draw (-2,0) node[below] {$\gamma_i'$};
\draw (2,0) node[below] {$\gamma_j'$};
\draw (0,0) node[below] {$\tau$};
\draw[red] (-4,0.1)--(-1,0.1);
\draw[blue] (4,0.1)--(1,0.1);
\draw[red] (-2.5,0) node[above] {$\gamma_i$};
\draw[blue] (2.5,0.01) node[above] {$\gamma_j$};

\draw (-4,0-2)--(4,0-2);
\draw (-4,0-2) node {$\bullet$};
\draw (4,0-2) node {$\bullet$};
\draw (-3,0-2) node {$\bullet$};
\draw (3,0-2) node {$\bullet$};
\draw (-1,0-2) node {$\bullet$};
\draw (1,0-2) node {$\bullet$};
\draw (-4,0-2) node[above] {$v_i'$};
\draw (4,0-2) node[above] {$v_j'$};
\draw (-1,0-2) node[above] {$v_i$};
\draw (1,0.01-2) node[above] {$v_j$};
\draw (-3.5,0-2) node[below] {$e_i'$};
\draw (3.5,0-2) node[below] {$e_j'$};
\draw (-2,0-2) node[below] {$\gamma_i'$};
\draw (2,0-2) node[above] {$\gamma_j'$};
\draw (0,-0.1-2) node[below] {$\tau$};
\draw[red] (-4,0.1-2)--(1,0.1-2);
\draw[blue] (4,-0.1-2)--(-1,-0.1-2);
\draw[red] (-2.5,0-2) node[above] {$\gamma_i$};
\draw[blue] (2.7,-0.1-2) node[below] {$\gamma_j$};
\end{tikzpicture}
\caption{The paths constructed in the proof of Lemma~\ref{Lem compact subset of boundary disjoint union}.}\label{Figure paths constructed}
\end{figure}
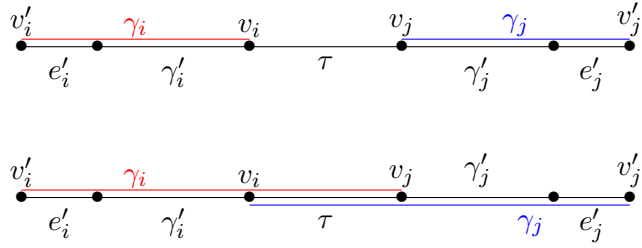

\dem See~Figure~\ref{Figure paths constructed} to follow the construction. Notice that we either have $\gamma_i \cap \gamma_j=\varnothing$ or $\gamma_i \cap \gamma_j \neq \varnothing$. In both cases, we construct a path $\tau$ and vertices $v_i,v_i',v_j,v_j'$ that we will use in the rest of the proof. First suppose that $\gamma_i \cap \gamma_j=\varnothing$. Let $\tau$ be the unoriented geodesic path in $T$ which realizes the distance between $\gamma_i$ and $\gamma_j$. Since, by assumption, $C(\gamma_i) \cap C(\gamma_j) \neq \varnothing$, the endpoints of $\tau$ are endpoints of $\gamma_i$ and $\gamma_j$. For every $p \in \{i,j\}$, let $v_p$ be the common endpoint of $\gamma_p$ and $\tau$ and let $v_p'$ be the other endpoint of $\gamma_p$. Suppose now that $\gamma_i \cap \gamma_j \neq \varnothing$. Then, since $C(\gamma_i) \cap C(\gamma_j) \neq \varnothing$ there exist three paths $\tau$, $a_i$ and $a_j$ such that, up to changing the orientation of $\gamma_i$ and $\gamma_j$, we have: $\gamma_i=a_i\tau$ and $\gamma_j=\tau a_j$. For every $p \in \{i,j\}$, let $v_p$ be the common endpoint of $a_p$ and $\tau$ and let $v_p'$ be the other endpoint of $a_p$. 

For every $p \in \{i,j\}$, let $e_p'$ be the edge of $\gamma_p$ adjacent to $v_p'$, which exists since $\gamma_p$ is not reduced to a vertex. For every $p \in \{i,j\}$, let $\gamma_p'$ be the edge path such that either $\gamma_p=\gamma_p'e_p'$ or $\gamma_p=e_p'\gamma_p'$. For every $p \in \{i,j\}$ and $\ell \in \{i,j\}-\{p\}$, let $\gamma_1^{(p)},\ldots,\gamma_{k_p}^{(p)}$ be the edge paths of $T$ which start at $v_p'$, which properly contain $\gamma_p$ and such that for every $s \in \{1,\ldots,k_p\}$, the endpoint of $\gamma_s^{(p)}$ distinct from $v_p'$ is at distance exactly $1$ from the minimal edge path of $T$ which contains $\tau$ and $\gamma_{\ell}'$. Note that for every $p \in \{i,j\}$ and $\ell \in \{i,j\}-\{p\}$, there exists a unique $s_p \in \{1,\ldots,k_p\}$ such that $\gamma_{s_p}^{(p)}$ contains $e_{\ell}'$. Note that for every $p \in \{i,j\}$, the integer $s_p$ is the unique integer $s \in \{1,\ldots,k_p\}$ such that $\gamma_s^{(p)}$ contains both $\gamma_i$ and $\gamma_j$. Note also that $\gamma_{s_i}^{(i)}=(\gamma_{s_j}^{(j)})^{-1}$.

We claim that the paths $\gamma_1^{(i)},\ldots,\gamma_{k_i}^{(i)},\gamma_1^{(j)},\ldots,\gamma_{k_j}^{(j)}$ satisfy the conclusion of the claim. Indeed, $(a)$ is satisfied by construction. We prove $(b)$. Let $p \in \{i,j\}$. Let $s,t \in \{1,\ldots,k_p\}$ be distinct. Then $\gamma_s^{(p)}$ and $\gamma_t^{(p)}$ share the path $\gamma_p$ as an initial segment. But, by construction of the paths $\gamma_s^{(p)}$ and $\gamma_t^{(p)}$, the endpoints of $\gamma_s^{(p)}$ and $\gamma_t^{(p)}$ distinct from $v_p'$ are at distance exactly $1$ from the minimal edge path of $T$ which contains $\tau$ and $\gamma_{\ell}'$. Therefore, the endpoint of $\gamma_s^{(p)}$ distinct from $v_p'$ is not contained in $\gamma_t^{(p)}$. Hence the subtree of $T$ generated by $\gamma_s^{(p)}$ and $\gamma_t^{(p)}$ is a tripod. This shows that $C(\gamma_s^{(p)}) \cap C(\gamma_t^{(p)})=\varnothing$ and this proves $(b)$. 

We now prove $(c)$. Let $s \in \{1,\ldots,k_i\}$ and let $t \in \{1,\ldots,k_j\}$. Suppose that we have $C(\gamma_s^{(i)}) \cap C(\gamma_t^{(i)}) \neq \varnothing$. Then there exists a path $\gamma'$ of $T$ such that $\gamma'$ contains both $\gamma_s^{(i)}$ and $\gamma_t^{(j)}$. Thus $\gamma'$ contains both $\gamma_i$ and $\gamma_j$. This implies that $\gamma_s^{(i)}=\gamma_{s_i}^{(i)}=(\gamma_{s_j}^{(j)})^{-1}=(\gamma_t^{(j)})^{-1}$ and that $C(\gamma_s^{(i)})=C(\gamma_t^{(j)})$. This proves $(c)$. Finally, the fact that $(d)$ holds follows from the fact that $C(\gamma)=\bigcup_{b \in ET, \gamma b \nsubseteq \gamma } C(\gamma b)$. This proves the claim.
\hfill\qedsymbol

\bigskip

For every $p \in \{i,j\}$, replace $\gamma_p$ by the paths $\gamma_1^{(p)},\ldots,\gamma_{k_i}^{(p)}$. Then we obtain a new set $\{\gamma_1',\ldots,\gamma_{\ell_1}'\}$ such that, by the point $(d)$ of the claim, $C=\cup_{i=1}^{\ell_1} C(\gamma_i')$. Recall that for every $p \in \{i,j\}$, we have $C(\gamma_p) \in \mathrm{Cyl}(\mathscr{C})$. By the point $(a)$ of the claim, for every $p \in \{i,j\}$ and every $s \in \{1,\ldots,k_i\}$, we have $\gamma_p \subseteq \gamma_s^{(p)}$. Therefore, we see that for every $p \in \{i,j\}$ and every $s \in \{1,\ldots,k_i\}$, we have $C(\gamma_s^{(p)}) \in \mathrm{Cyl}(\mathscr{C})$. Hence the set $\{\gamma_1',\ldots,\gamma_{\ell_1}'\}$ satisfies $(1)$. Point $(a)$ of the claim also implies that, for every $m' \in \{1,\ldots,\ell\}$, and every $p \in \{i,j\}$, if $C(\gamma_{m'}) \cap C(\gamma_p)=\varnothing$ then for every $s \in \{1,\ldots,k_p\}$, we have $C(\gamma_{m'}) \cap C(\gamma_s^{(p)})=\varnothing$. Combined with points $(b)$ and $(c)$ of the claim, we see that the number of distinct elements $m_1,m_2 \in \{1,\ldots,\ell_1\}$ such that $C(\gamma_{m_1}) \cap C(\gamma_{m_2}) \neq \varnothing$ is strictly less than $m$. An inductive argument then concludes the proof.
\hfill\qedsymbol

\bigskip

We denote by $F_n-\mathcal{A}$ the subset of $F_n$ consisting in every element $w \in F_n$ such that $C(\gamma_{w}) \in \mathrm{Cyl}(\mathscr{C})$. Note that $F_n-\mathcal{A}$ is closed under inversion since $\mathscr{C}$ is closed under inversion by Lemma~\ref{Lem finite set of words determines vertex group system}. The next lemma gives a criterion to extend some functions defined on $F_n-\mathcal{A}$ to a relative current in $\Curr(F_n,\mathcal{A})$ (see~\cite[Lemma~3.9]{gupta2017relative} for the free factor system case). First we need some definitions. 

Let $w \in F_n$, and let $k \in \NN^*$. A \emph{length $k$ extension of $w$} is a word $w'=wx_1\ldots x_k$ where for every $i \in \{1,\ldots,k-1\}$, we have $x_i \neq x_{i+1}^{-1}$ and $x_1$ is not the inverse of the last letter of $w$. An \emph{extension of $w$} is a word $w'$ such that there exists $k \in \NN^*$ such that $w'$ is a length $k$ extension of $w$.

\begin{lem}\label{Lem Kolmogorov extension}
Let $\eta \colon F_n -\mathcal{A} \to \RR_+$ be a function invariant under inversion and which satisfies, for every $w \in F_n -\mathcal{A}$:
\begin{equation}\label{Equation Kolmogorov}
\eta(w)=\sum\limits_{v \text{ is a length one extension of }w}\eta(v).
\end{equation}
There exists a unique element $\widetilde{\eta} \in \Curr(F_n,\mathcal{A})$ such that for every element $w \in F_n-\mathcal{A}$, we have $$\eta(w)=\widetilde{\eta}(C(\gamma_w)).$$
\end{lem}

\dem Since $\partial^2(F_n,\mathcal{A})$ is totally disconnected and locally compact by Lemma~\ref{Lem boundary loc compact}, and since a relative current is a Radon measure, a relative current is uniquely determined by its values on compact open subsets of $\partial^2(F_n,\mathcal{A})$. Let $C$ be a compact open subset of $\partial^2(F_n,\mathcal{A})$. By Lemma~\ref{Lem compact subset of boundary disjoint union}, the subset $C$ is a disjoint union of cylinders of finitely many geodesic edge paths $\gamma_1,\ldots,\gamma_k$ such that for every $i \in \{1,\ldots,k\}$, we have $C(\gamma_i) \in \mathrm{Cyl}(\mathscr{C})$. For every $i \in \{1,\ldots,k\}$, let $g_i$ be the element of $F_n$ which is the label of $\gamma_i$. For every $i \in \{1,\ldots,r\}$, since $g_i$ contains an element of $\mathscr{C}$ as a subword, we have $g_i \in F_n-\mathcal{A}$. Hence we can set $\widetilde{\eta}(C)=\sum_{i=1}^k \eta(g_i)$. We claim that the value $\widetilde{\eta}(C)$ does not depend on the choice of the paths $\gamma_i$. Indeed, let $\alpha_1,\ldots,\alpha_{\ell}$ be another set of geodesic edge paths given by Lemma~\ref{Lem compact subset of boundary disjoint union} and let $h_1,\ldots,h_{\ell}$ be the corresponding elements in $F_n$. Note that for every $i \in \{1,\ldots,k\}$ and every $j \in \{1,\ldots,\ell\}$ such that $C(\gamma_i) \cap C(\alpha_j) \neq \varnothing$, we have $C(\gamma_i) \cap C(\alpha_j)=C(\beta_{i,j})$, where $\beta_{i,j}$ is a minimal edge path in $T$ that contains both $\gamma_i$ and $\alpha_j$. 

We claim that for every $i \in \{1,\ldots,k\}$, there do not exist distinct $j_1,j_2 \in \{1,\ldots,\ell\}$ and paths $a_1$ and $a_2$ such that $\beta_{i,j_1}=a_1\gamma_i$ and $\beta_{i,j_2}=\gamma_ia_2$. Indeed, otherwise the path $a_1\gamma_ia_2$ is a finite path that contains both $\alpha_{j_1}$ and $\alpha_{j_2}$. Hence $C(\alpha_{j_1}) \cap C(\alpha_{j_2}) \neq \varnothing$, a contradiction. The claim follows.

For every $i \in \{1,\ldots,k\}$ and every $j \in \{1,\ldots,\ell\}$ such that $C(\gamma_i) \cap C(\alpha_j) \neq \varnothing$, let $g_{i,j}$ be an element in $F_n$ corresponding to $\beta_{i,j}$. By the above claim, for every $i \in \{1,\ldots,k\}$, one of the following holds:

\medskip

\noindent{$(a)$ } for every $j \in \{1,\ldots,\ell\}$ such that $C(\gamma_i) \cap C(\alpha_j) \neq \varnothing$, the element $g_{i,j}$ is an extension of $g_i$;

\medskip

\noindent{$(b)$ } for every $j \in \{1,\ldots,\ell\}$ such that $C(\gamma_i) \cap C(\alpha_j) \neq \varnothing$, the element $g_{i,j}^{-1}$ is an extension of $g_i^{-1}$.

Since $\eta$ is invariant under inversion, we may suppose that for every $i \in \{1,\ldots,k\}$, and for every $j \in \{1,\ldots,\ell\}$ such that $C(\gamma_i) \cap C(\alpha_j) \neq \varnothing$, the element $g_{i,j}$ is an extension of $g_i$. Thus for every $j \in \{1,\ldots,\ell\}$, and for every $i \in \{1,\ldots,k\}$ such that $C(\gamma_i) \cap C(\alpha_j) \neq \varnothing$, the element $g_{i,j}^{-1}$ is an extension of $h_j^{-1}$.

Note that, since $C= \cup_{i=1}^k C(\gamma_i)=\cup_{j=1}^{\ell} C(\alpha_j)$, for every $i \in \{1,\ldots,k\}$, the subset $C(\gamma_i)$ is covered by a disjoint union of finitely many $C(\alpha_j)$. Hence, for every $i \in \{1,\ldots,k\}$, Equation~\eqref{Equation Kolmogorov} implies that: $$\eta(g_i)=\sum\limits_{j\;|\; C(\gamma_i)\cap C(\alpha_j) \neq \varnothing }\eta(g_{i,j}). $$

Similarly, for every $j \in \{1,\ldots,\ell\}$, we have:
$$\eta(h_j^{-1})=\sum\limits_{i\;|\; C(\gamma_i)\cap C(\alpha_j) \neq \varnothing }\eta(g_{i,j}^{-1}). $$

Thus, since $\eta$ is invariant under inversion, we have:
$$
\sum_{j=1}^{\ell} \eta(h_j)=\sum_{j=1}^{\ell} \eta(h_j^{-1})=\sum_{j=1}^{\ell}\sum_{i\;|\; C(\gamma_i)\cap C(\alpha_j) \neq \varnothing }\eta(g_{i,j}^{-1})=\sum_{i=1}^{k}\sum_{j\;|\; C(\gamma_i)\cap C(\alpha_j) \neq \varnothing }\eta(g_{i,j})=\sum_{i=1}^k \eta(g_i).
$$

Hence the value of $\widetilde{\eta}(C)$ does not depend on the choice of the paths $\gamma_i$.

Therefore $\widetilde{\eta}$ is an additive, $F_n$-invariant and nonnegative function on the set of compact open subsets of $\partial^2(F_n,\mathcal{A})$. We claim that $\widetilde{\eta}$ is in fact $\sigma$-additive. Indeed, by \cite[Proposition~1.2.6]{Cohn80}, it suffices to prove that for every decreasing sequence $(C_n)_{n \in \NN}$ of compact open subsets of $\partial^2(F_n,\mathcal{A})$ such that $\bigcap_{n \in \NN} C_n=\varnothing$, we have $\lim_{n \to \infty} \widetilde{\eta}(C_n)=0$. But since a decreasing sequence of nonempty compact subsets is a nonempty compact subset, there exists $n \in \NN$ such that $C_n=\varnothing$. This proves the claim. By Carathéodory extension theorem (see~\cite[Proposition~1.2.6, Theorem~1.3.6]{Cohn80}), the function $\widetilde{\eta}$ has a unique extension as a Radon measure on the $\sigma$-algebra of Borel sets of $\partial^2(F_n,\mathcal{A})$. 
\hfill\qedsymbol

\bigskip

Let $$\PCurr(F_n,\mathcal{A})=\left(\Curr(F_n,\mathcal{A})-\{0\}\right)/\RR_+^*$$ be the set of projectivized relative currents (where $\RR_+^*$ acts on $\Curr(F_n,\mathcal{A})$ by homothety), equipped with the quotient topology which is metrizable. The next result is a generalization of \cite[Lemma~3.11]{gupta2017relative}.

\begin{lem}\label{Lem PCurr compact}
The metrisable space $\PCurr(F_n,\mathcal{A})$ is compact.
\end{lem}

\dem Let $([\eta_n])_{n \in \NN}$ be a sequence of projective currents relative to $\mathcal{A}$. We prove that it has a convergent subsequence. Let $\mathscr{C}$ be the finite set given by Lemma~\ref{Lem finite set of words determines vertex group system}. For every $n \in \NN$, let $\eta_n$ be a representative of $[\eta_n]$ such that, for every $w \in \mathscr{C}$, we have $\eta(C(\gamma_w)) \leq 1$, with equality for some $w \in \mathscr{C}$, independent of $n$ up to extraction. The set $\mathscr{C}$ being finite, there exists a subsequence $(n_k)_{k \in \NN}$ such that for every $u \in \mathscr{C}$, the sequence $(\eta_{n_k}(C(\gamma_u)))_{k \in \NN}$ converges. Moreover, there exists $u_0 \in \mathscr{C}$ such that the limit $\lim_{k \to \infty}(\eta_{n_k}(C(\gamma_{u_0})))_{k \in \NN}$ is not equal to zero. Let $w \in F_n$ be such that $C(\gamma_w) \in \mathrm{Cyl}(\mathscr{C})$. There exists $u_w \in \mathscr{C}$ such that $u_w$ is a subword of $w$. Therefore, for every $k \in \NN$, we have $$\eta_{n_k}(C(\gamma_w)) \leq \eta_{n_k}(C(\gamma_{u_w})) \leq 1.$$ Therefore, for every element $w \in F_n-\mathcal{A}$, the sequence $(\eta_{n_k}C((\gamma_w)))_{k \in \NN}$ has a convergent subsequence. By a diagonal argument, up to extraction, for every $C(\gamma_w) \in \mathrm{Cyl}(\mathscr{C})$, the sequence $(\eta_{n_k}(C(\gamma_w)))_{k \in \NN}$ converges. Moreover, there exists $C(\gamma_w) \in \mathrm{Cyl}(\mathscr{C})$ such that $(\eta_{n_k}(C(\gamma_w)))_{k \in \NN}$ converges to a nonzero element. 

Let $\eta \colon F_n -\mathcal{A} \to \RR_+$ be the function defined by, for every $w \in F_n-\mathcal{A}$: $$\eta(w)=\lim_{k \to \infty} \eta_{n_k}(C(\gamma_w)).$$

Since for every $k \in \NN$, the function $\eta_{n_k}$ is a relative current, the function $\eta$ satisfies the assumptions of Lemma~\ref{Lem Kolmogorov extension}. Therefore, by Lemma~\ref{Lem Kolmogorov extension}, there exists a unique relative current $\widetilde{\eta} \in \Curr(F_n,\mathcal{A})$ such that for every element $w \in F_n-\mathcal{A}$, we have $$\eta(w)=\widetilde{\eta}(C(\gamma_w)).$$ Hence $([\eta_{n_k}])_{k \in \NN}$ converges to $[\widetilde{\eta}]$.
\hfill\qedsymbol

\section{Density of rational currents}

In this section, let $n \geq 3$. Let $r \in \NN$ and let $\mathcal{A}=\{[A_1],\ldots,[A_r]\}$ be a malnormal subgroup system of $F_n$. Let $\mathcal{B}$, $T$, $L \in \NN^*$, $\Gamma_1,\ldots,\Gamma_r$, $\mathscr{C}$ be as above Lemma~\ref{Lem finite set of words determines vertex group system}. Let $\ell \colon F_n \to \NN$ be the length function corresponding to $\mathcal{B}$. 

Every conjugacy class of nonperipheral element $g \in F_n$ determines a relative current $\eta_g$ as follows. Suppose first that $g$ is \emph{root-free}, that is, $g$ is not a proper power of any element in $F_n$. Let $\gamma$ be a finite geodesic path in the Cayley graph $T$ such that $C(\gamma) \in \mathrm{Cyl}(\mathscr{C})$. Then $\eta_g(C(\gamma))$ is the number of unoriented translation axes in $T$ of conjugates of $g$ that contain the path $\gamma$. If $g=h^k$ with $k \geq 2$ and $h$ root-free, we set $\eta_g=k \; \eta_h$. Such currents are called \emph{rational currents}. Note that for every nonperipheral element $g \in F_n$, the current $\eta_g$ only depends on the conjugacy class of $g$. Therefore, we can talk about rational currents induced by conjugacy classes of nonperipheral elements of $F_n$ and write $\eta_{[g]}$ for the rational current associated with the conjugacy class of a nonperipheral element $g \in F_n$. We prove the following proposition.

\begin{prop}\label{Prop density rational currents}
Let $n \geq 3$ and let $\mathcal{A}$ be a malnormal subgroup system of $F_n$. The set of projectivized rational currents induced by conjugacy classes of nonperipheral elements of $F_n$ is dense in $\PCurr(F_n,\mathcal{A})$.
\end{prop}

We follow Gupta's proof (\cite[Proposition~3.12]{gupta2017relative}) in the special case of free factor systems. The proof consists in approximating currents in 
$\PCurr(F_n,\mathcal{A})$ with \emph{signed measured currents} on $\partial^2F_n$, which are $F_n$-invariant and $\sigma$-additive real-valued functions on the set of Borel subsets of $\partial^2F_n$. We will then conclude using the following lemma, due to Martin (see also~\cite[Lemma~3.15]{gupta2017relative}).

\begin{lem}\cite[Lemma~15]{Martin95}\label{Lem approximation signed measured currents by rational currents}
Let $n \geq 3$. Suppose that $\mathcal{A}=\varnothing$. Let $k' \geq 1$, let $k \geq 2$ with $k' \leq k$ and let $\eta$ be a signed measured current such that, for every $w \in F_n$ with $k' \leq \ell(w) \leq k$, we have $\eta(C(\gamma_w)) \geq 0$. Let $P=2n(2n-1)^{2n(2n-1)^{k-2}}$. If there exists $w_0 \in F_n$ such that $\ell(w_0)=k$ and $\eta(C(\gamma_{w_0})) \geq P$, then there exists $\alpha \in F_n-\{e\}$ such that, for every $w \in F_n$ with $k' \leq \ell(w) \leq k$, we have $\eta(C(\gamma_w)) \geq \eta_{[\alpha]}(C(\gamma_w))$.
\end{lem}

\begin{rmq}\label{Rmq lem approximation by rational currents}
\noindent{$(1)$ } The hypotheses in \cite[Lemma~15]{Martin95} requires that $k'=1$. However, the proof of Martin works by studying words of length exactly $k$ and then extend the result to words of length at most $k$ by additivity of the measures. Thus the proof with $k'>1$ is identical. 

\medskip

\noindent{$(2)$ } For the rational current $\eta_{[\alpha]}$ constructed in Lemma~\ref{Lem approximation signed measured currents by rational currents}, there exists  $w \in F_n$ with $k' \leq \ell(w) \leq k$ such that $\eta_{[\alpha]}(C(\gamma_w))>0$.
\end{rmq}

Recall that $\mathrm{Cyl}(\mathscr{C})$ is the set of cylinder sets of the form $C(\gamma_w)$, where $w$ is a word of $F_n$ containing a word of $\mathscr{C}$ as a subword. Let $\eta_0 \in \Curr(F_n,\mathcal{A})$ and let $k \geq L+2$. Let $\eta$ be a signed measured current such that, for every element $w \in F_n$ with $C(\gamma_w) \in \mathrm{Cyl}(\mathscr{C})$, we have $\eta(C(\gamma_w))=\eta_0(C(\gamma_w))$ and for every element $w \in F_n$ of length between $L+2$ and $k$, we have $\eta(C(\gamma_w)) \geq 0$. Then $\eta$ is called a \emph{$k$-extension of $\eta_0$}. The key lemma in order to prove Proposition~\ref{Prop density rational currents} is the following result (see \cite[Lemma~3.15]{gupta2017relative} for the same statement in the particular case free factor systems):

\begin{lem}\label{Lem k extension}
Let $\eta_0$ be a relative current and let $k \geq L+2$. There exists a signed measured current $\eta \colon \partial^2 F_n \to \RR$ which is a $k$-extension of $\eta_0$.
\end{lem}

Let $\eta_0$ be a relative current. In order to prove Lemma~\ref{Lem k extension}, we need some preliminary results. We follow \cite[Section~8.1]{gupta2017relative}. For $k \in \NN^*$, let $S_k$ be the set of elements of $F_n$ of length $k$ which do not contain an element of $\mathscr{C}$ as a subword. Note that, since $\mathscr{C}$ is closed under inversion by Lemma~\ref{Lem finite set of words determines vertex group system}, we see that, for every $k \in \NN^*$, the set $S_k$ is closed under inversion. For $k=0$, we set $S_0=\{e\}$. Note also that, if $k<L+2$, then $S_k$ contains all words of length $k$ since every element of $\mathscr{C}$ has length equal to $L+2$.

\medskip

\begin{lem}\label{Lem sets S_k} $(1)$ If $\mathcal{A} \neq \varnothing$, for every $k \in \NN^*$, the set $S_k$ is not empty. 

\medskip

\noindent{$(2)$ } For every $k \geq L+2$ and every $w \in S_k$, there exist $w' \in S_{k+1}$, $i \in \{1,\ldots,r\}$, $g \in F_n$ and $a \in gA_ig^{-1}$ such that $w'$ is a length $1$ extension of $w$ and $a$ is an extension of $w'$.

\end{lem}

\dem $(1)$ Since the group $A_1$ is infinite, the corresponding minimal subtree $T_{A_1}$ is infinite. Recall that the tree $T_{A_1}$ is supposed to contain the origin $e$ of $T$. Let $\gamma$ be a geodesic path contained in $T_{A_1}$, starting from $e$ and of length equal to $k$, and let $h \in F_n$ be the corresponding element of $F_n$. Then there exists $a \in A_1$ such that $a$ is an extension of $h$. We have $h \in S_k$ as otherwise $a$ would contradict Lemma~\ref{Lem finite set of words determines vertex group system}~$(3)$. This proves $(1)$.

\medskip

\noindent{$(2)$ } Let $k \geq L+2$ and let $w \in S_k$. By Lemma~\ref{Lem finite set of words determines vertex group system}~$(1)$, there exist $i \in \{1,\ldots,r\}$ and $g \in F_n$ such that $\gamma_w$ is contained in $T_{gA_ig^{-1}}$. As $T_{gA_ig^{-1}}$ does not contain any univalent vertex, there exists a geodesic ray $c$ in $T_{gA_ig^{-1}}$ starting from $e$ which contains the path $\gamma_w$. Let $\gamma'$ be the geodesic path in $c$ of length $k+1$ containing $\gamma_w$, and let $w'$ be the corresponding element in $F_n$. Then $w' \in S_{k+1}$ and $w'$ is a length $1$ extension of $w$. This proves $(2)$ and this concludes the proof.
\hfill\qedsymbol

\bigskip

Let $k \geq L+2$. Let $S_k^0$ be a subset of $S_k$ (chosen once and for all) such that for every $w \in S_k$ exactly one of $w$ or $w^{-1}$ appears in $S_k^0$. In what follows, we adopt the convention that whenever an extension of a word $w$ by a letter $b \in \mathcal{B}$ is written as $wb$ (resp. $bw$), we assume that $b$ is not the inverse of the last letter (resp. first letter) of the word $w$.

In order to construct the signed measured current which satisfies the conclusion of Lemma~\ref{Lem k extension}, we will define a signed measured current on cylinders of words in $S_{k-1}$ and use those values together with the additivity laws in order to define $\eta$ on cylinders of words of length $k$. First we set $\eta(C(\gamma_b))=1$ for every letter $b$ of $\mathcal{B}$ not contained in $\mathscr{C}$. By induction, assume that for every element $v \in S_{k-1}$, the value $\eta(C(\gamma_v))$ is defined. By additivity of a signed measured current, for every $v \in S_{k-1}^0$, we want to have:

$$\begin{array}{c}
\eta(C(\gamma_v))=\sum\limits_{b \in \mathcal{B},vb \in S_k} \eta(C(\gamma_{vb}))+\sum\limits_{b \in \mathcal{B},vb \notin S_k} \eta_0(C(\gamma_{vb})) \\
\eta(C(\gamma_{v^{-1}}))=\sum\limits_{b \in \mathcal{B},v^{-1}b \in S_k} \eta(C(\gamma_{v^{-1}b}))+\sum\limits_{b \in \mathcal{B},v^{-1}b \notin S_k} \eta_0(C(\gamma_{v^{-1}b}))
\end{array}
$$

Since $\eta$ is invariant under taking inverses, the equation obtained by using forward extensions of $v^{-1}$ is the same one as the equation obtained by using backward extensions of $v$. After rearranging the equations in order to have the unknown terms on the left hand side, we obtain:

\begin{equation}\label{Equation E_{k-1}}
\begin{array}{c}
\sum\limits_{b \in \mathcal{B},vb \in S_k} \eta(C(\gamma_{vb}))=\sum\limits_{b \in \mathcal{B},vb \notin S_k} \eta_0(C(\gamma_{vb})) -\eta(C(\gamma_v))=c_v\\
\sum\limits_{b \in \mathcal{B},v^{-1}b \in S_k} \eta(C(\gamma_{v^{-1}b}))=\sum\limits_{b \in \mathcal{B},v^{-1}b \notin S_k} \eta_0(C(\gamma_{v^{-1}b}))-\eta(C(\gamma_{v^{-1}}))=c_{v^{-1}}.
\end{array}
\end{equation}

Since $\eta$ is invariant under taking inverse, this shows that there are $|S_{k-1}|$ equations in $|S_k|/2=|S_k^0|$ variables.

Denote the system of equations~\eqref{Equation E_{k-1}} by $E_{k-1}^1$. These are equations obtained from length $1$ extensions of words in $S_{k-1}$. Similarly, for every $i \in \{1,\ldots,k-1\}$, we define $E_{k-i}^i$ as the system of equations obtained from length $i$ extensions of words in $S_{k-i}$.

Let $[M|c]$ be the augmented matrix for the system of equations $E_{k-1}^1$ with rows labeled by words in $S_{k-1}$, columns by words in $S_k^0$ and such that for every $w \in S_k^0$ and every $v \in S_{k-1}$, we have $M_{v,w}=1$ if there exists $b \in \mathcal{B}$ such that $w=vb$ or $w^{-1}=vb$; and $M_{v,w}=0$ otherwise. Let $c$ be the column vector indexed by words in $S_{k-1}$ such that for every $v \in S_{k-1}$, the coordinate of $c$ at $v$ is equal to $c_v$. If $v \in S_{k-1}$, we will denote by $r_v$ the corresponding row vector of $M$. Observe that each column has exactly two entries which are equal to $1$. Indeed, $M_{v,w}$ is equal to $1$ exactly when $w$ or $w^{-1}$ is a length $1$ extension of $v$. Observe also that any two distinct row vectors $r_{v_1}$ and $r_{v_2}$ can have at most one common coordinate which is equal to $1$. Indeed, let $w \in S_k^0$ be such that $M_{v_1,w}=M_{v_2,w}=1$. Then there exist $b_1,b_2 \in \mathcal{B}$ such that $w=v_1b_1$ or $w=b_1^{-1}v_1^{-1}$ and $w=v_2b_2$ or $w=b_2^{-1}v_2^{-1}$. Therefore, the word $v_1$ starts with $b_2^{-1}$ and $v_2$ starts with $b_1^{-1}$. This shows that $w$ is uniquely determined. 

The next lemma is the same one as \cite[Lemma~8.2]{gupta2017relative} in the special case of free factor systems.

\begin{lem}\label{Lem equations k extension}
\noindent{$(1)$ } For every $i \geq 1$, an equation in the system $E_{k-i-1}^{i+1}$ is a linear combination of equations in the system $E_{k-i}^i$. Thus it is sufficient to look at the system $E_{k-1}^1$ in order to obtain every constraint satisfied by $\eta(C(\gamma_w))$ for every $w \in S_k^0$.

\medskip

\noindent{$(2)$ } Let $u \in S_{k-2}$. Then the following two linear combinations of rows of $M$ are equal: 
\begin{equation}\label{Equation lem k extension 2}
\sum\limits_{b \in \mathcal{B}, bu \in S_{k-1}}r_{bu}=\sum\limits_{b \in \mathcal{B}, bu^{-1} \in S_{k-1}}r_{bu^{-1}}.
\end{equation}

\medskip

\noindent{$(3)$ } Every relation among the rows of $M$ is a linear combination of relations in the set of relations~\eqref{Equation lem k extension 2} where $u$ varies in $S_{k-2}$.

\medskip

\noindent{$(4)$ } We have $$\sum\limits_{b \in \mathcal{B}, bu \in S_{k-1}}c_{bu}=\sum\limits_{b \in \mathcal{B}, bu^{-1} \in S_{k-1}}c_{bu^{-1}},$$ where for every $v \in S_{k-1}$, $c_v$ is given by Equation~\eqref{Equation E_{k-1}}.

\medskip

\noindent{$(5)$ } The system of equations $E_{k-1}^1$ is consistent and hence has a solution. Thus we can define $\eta$ on words of length $k$.
\end{lem}

\dem \noindent{$(1)$ } Let $i \geq 1$ and $u \in S_{k-i-1}$. Then by the system $E_{k-i-1}^1$ $$ \eta(C(\gamma_u))=\sum\limits_{b \in \mathcal{B}} \eta(C(\gamma_{ub})).$$ By the equations in $E_{k-i}^{i}$, we have, for every $b \in \mathcal{B}$: $$\eta(C(\gamma_{ub}))=\sum\limits_{y \in F_n, \ell(y)=i}\eta(C(\gamma_{uby})).$$ Adding all these equations over $b \in \mathcal{B}$, we have:
$$\eta(C(\gamma_u))=\sum\limits_{b,y \in F_n, \ell(b)=1, \ell(y)=i} \eta(C(\gamma_{uby}))=\sum\limits_{z \in F_n, \ell(z)=i+1} \eta(C(\gamma_{uz})).$$ Thus we have recovered an equation in $E_{k-i-1}^{i+1}$ as a linear combination of equations in $E_{k-i}^i$. This proves $(1)$.

\medskip

\noindent{$(2)$ } Let $u \in S_{k-2}$ and let $w \in S_k^0$. For every $b \in \mathcal{B}$ such that $bu \in S_{k-1}$, we have $M_{bu,w} \neq 0$ exactly when there exists $y \in \mathcal{B}$ such that $w=buy^{-1}$ or $w=yu^{-1}b^{-1}$ (recall that the basis $\mathcal{B}$ is supposed to be symmetric). Therefore, if $M_{bu,w} \neq 0$, there exists a unique $y \in \mathcal{B}$ such that $M_{yu^{-1},w} \neq 0$. This proves $(2)$.

\medskip

\noindent{$(3)$ } Let $R$ be a relation given by $\sum_{v \in S_{k-1}} d_vr_v=0$, where $d_v \in \RR$. Suppose that the number of terms in the sum associated with $R$ is minimal. Such an assumption is possible as every relation is a linear combination of relations whose number of terms is minimal. We can rescale the equation so that there exist $b \in \mathcal{B}$ and $u \in S_{k-2}$ such that $d_{bu}=1$. For every $y \in \mathcal{B}$ such that $buy^{-1} \in S_{k}^0$, we have $$M_{bu,buy^{-1}}=M_{yu^{-1},buy^{-1}}=1.$$ This implies, as explained above the lemma, that the rows $r_{bu}$ and $r_{yu^{-1}}$ share exactly one common nonzero coordinate, which is $buy^{-1}$. Moreover, the rows $r_{bu}$ and $r_{yu^{-1}}$ are the only rows which have a nonzero coordinate in $buy^{-1}$. This shows that $d_{yu^{-1}}=-1$.

Let $y \in \mathcal{B}$ be such that $yu^{-1} \in S_{k-1}$. For every $z \in \mathcal{B}$ such that $yu^{-1}z \in S_k^0$, we have $M_{yu^{-1},yu^{-1}z}=M_{z^{-1}u,yu^{-1}z}=1$. Thus we have $d_{z^{-1}u}=1$.  Therefore we see that $$\sum\limits_{b \in \mathcal{B},bu \in S_{k-1}} d_{bu}r_{bu}- \sum\limits_{y \in \mathcal{B},yu^{-1} \in S_{k-1}} d_{yu^{-1}}r_{yu^{-1}}=\sum\limits_{b \in \mathcal{B},bu \in S_{k-1}} r_{bu}- \sum\limits_{y \in \mathcal{B},yu^{-1} \in S_{k-1}} r_{yu^{-1}}=0.$$ Hence the minimal relation $R$ is just $$\sum\limits_{b \in \mathcal{B},bu \in S_{k-1}} r_{bu}- \sum\limits_{y \in \mathcal{B},yu^{-1} \in S_{k-1}} r_{yu^{-1}}=0.$$

\medskip

\noindent{$(4)$ } Let $u \in S_{k-2}$. We have, by the definition of $c_v$:

$$\begin{array}{ccl}
-\sum\limits_{b \in \mathcal{B}, bu \in S_{k-1}} c_{bu}&=& \sum\limits_{b \in \mathcal{B}, bu \in S_{k-1}} \eta(C(\gamma_{bu})) - \sum\limits_{b,y \in \mathcal{B}, bu \in S_{k-1},buy \notin S_{k}}\eta(C(\gamma_{buy})) \\
{} & = & \eta(C(\gamma_u))-\sum\limits_{b \in \mathcal{B}, bu \notin S_{k-1}} \eta(C(\gamma_{bu}))-\sum\limits_{b,y \in \mathcal{B}, bu \in S_{k-1},buy \notin S_{k}}\eta(C(\gamma_{buy})) \\
{} & = & \eta(C(\gamma_u))-\sum\limits_{b,y \in \mathcal{B}, bu \notin S_{k-1}} \eta(C(\gamma_{buy}))-\sum\limits_{b,y \in \mathcal{B}, bu \in S_{k-1},buy \notin S_{k}}\eta(C(\gamma_{buy})).
\end{array}
$$

Note that we have: 

\begin{equation}\label{Equation noninverse}
\sum\limits_{b,y \in \mathcal{B}, bu \notin S_{k-1}} \eta(C(\gamma_{buy}))=\sum\limits_{b,y \in \mathcal{B}, bu \notin S_{k-1},uy \in S_{k-1}} \eta(C(\gamma_{buy}))+\sum\limits_{b,y \in \mathcal{B}, bu \notin S_{k-1},uy \notin S_{k-1}} \eta(C(\gamma_{buy}))
\end{equation}

Similarly, we have:

$$
-\sum\limits_{b \in \mathcal{B}, bu^{-1} \in S_{k-1}} c_{bu^{-1}}= \eta(C(\gamma_u))-\sum\limits_{b,y \in \mathcal{B}, bu^{-1} \notin S_{k-1}} \eta(C(\gamma_{bu^{-1}y}))-\sum\limits_{b,y \in \mathcal{B}, bu^{-1} \in S_{k-1},bu^{-1}y \notin S_{k}}\eta(C(\gamma_{bu^{-1}y})).$$

The right hand side is also equal to:

$$ \eta(C(\gamma_u))-\sum\limits_{b,y \in \mathcal{B}, ub^{-1} \notin S_{k-1}} \eta(C(\gamma_{y^{-1}ub^{-1}}))-\sum\limits_{b,y \in \mathcal{B}, ub^{-1} \in S_{k-1},y^{-1}ub^{-1} \notin S_{k}}\eta(C(\gamma_{y^{-1}ub^{-1}})).
$$

Observe that the sum $\sum\limits_{b,y \in \mathcal{B}, ub^{-1} \notin S_{k-1}} \eta(C(\gamma_{y^{-1}ub^{-1}}))$ equals: \begin{equation}\label{Equation inverse}
\sum\limits_{b,y \in \mathcal{B}, ub^{-1} \notin S_{k-1},y^{-1}u \in S_{k-1}} \eta(C(\gamma_{y^{-1}ub^{-1}})) + \sum\limits_{b,y \in \mathcal{B}, ub^{-1} \notin S_{k-1},y^{-1}u \notin S_{k-1}} \eta(C(\gamma_{y^{-1}ub^{-1}})).
\end{equation}

Suppose first that $k \leq L+2$. Then $S_{k-1}$ contains all words of length $k-1$. Hence we have $$-\sum\limits_{b \in \mathcal{B}, bu \in S_{k-1}} c_{bu}=\eta(C(\gamma_u))-\sum\limits_{b,y \in \mathcal{B},buy \notin S_{k}}\eta(C(\gamma_{buy}))$$ and $$-\sum\limits_{b \in \mathcal{B}, bu^{-1} \in S_{k-1}} c_{bu^{-1}}=\eta(C(\gamma_u))-\sum\limits_{b,y \in \mathcal{B},y^{-1}ub^{-1} \notin S_{k}}\eta(C(\gamma_{y^{-1}ub^{-1}})),$$ so that Assertion~$(4)$ holds in this case with $y=b^{-1}$.

Suppose now that $k>L+2$. Then since every element of $\mathscr{C}$ has length equal to $L+2$, an element of $\mathscr{C}$ contained in a word $x$ of length $k$ is properly contained in $x$. Hence if $b,y \in \mathcal{B}$ are such that $bu \in S_{k-1}$ and $buy \notin S_k$, then $uy \notin S_{k-1}$. Thus, we see that: \begin{equation}\label{Equation new 1}
\sum\limits_{b,y \in \mathcal{B}, bu \in S_{k-1},buy \notin S_{k}}\eta(C(\gamma_{buy}))=\sum\limits_{b,y \in \mathcal{B}, bu \in S_{k-1},uy \notin S_{k-1}}\eta(C(\gamma_{buy})).
\end{equation} Similarly, we have: 

\begin{equation}\label{Equation new 2}
\sum\limits_{b,y \in \mathcal{B}, ub^{-1} \in S_{k-1},y^{-1}ub^{-1} \notin S_{k}}\eta(C(\gamma_{y^{-1}ub^{-1}}))=\sum\limits_{b,y \in \mathcal{B}, ub^{-1} \in S_{k-1},y^{-1}u \notin S_{k-1}}\eta(C(\gamma_{y^{-1}ub^{-1}})).
\end{equation}

Using Equations~\eqref{Equation noninverse},~\eqref{Equation inverse},~\eqref{Equation new 1} and~\eqref{Equation new 2} with $y=b^{-1}$, we see that 

$$
\begin{array}{c}
\sum\limits_{b,y \in \mathcal{B}, bu \notin S_{k-1}} \eta(C(\gamma_{buy}))+\sum\limits_{b,y \in \mathcal{B}, bu \in S_{k-1},buy \notin S_{k}}\eta(C(\gamma_{buy})) \\
=\sum\limits_{b,y \in \mathcal{B}, ub^{-1} \notin S_{k-1}} \eta(C(\gamma_{y^{-1}ub^{-1}}))+\sum\limits_{b,y \in \mathcal{B}, ub^{-1} \in S_{k-1},y^{-1}ub^{-1} \notin S_{k}}\eta(C(\gamma_{y^{-1}ub^{-1}})).
\end{array}
$$

This shows that $$\sum\limits_{b \in \mathcal{B}, bu \in S_{k-1}}c_{bu}=\sum\limits_{b \in \mathcal{B}, bu^{-1} \in S_{k-1}}c_{bu^{-1}},$$ and this proves~$(4)$.

\bigskip

\noindent{$(5)$ } By Assertions~$(3)$ and $(4)$, if $R$ is a linear combination of relations among the rows of $M$ equal to zero, then the corresponding linear combination among coordinates of the vector $c$ is also equal to zero. Therefore, the system $[M|c]$ has a solution.
\hfill\qedsymbol

\bigskip

\noindent{\it Proof of Lemma~\ref{Lem k extension} } Let $\eta_0$ be a relative current. By Lemma~\ref{Lem equations k extension}, there exists a signed measured current $\eta$ such that, for every element $w$ of $F_n$ which satisfies $C(\gamma_w) \in \mathrm{Cyl}(\mathscr{C})$, we have $\eta_0(C(\gamma_w))=\eta(C(\gamma_w))$. This extension is not necessarily nonnegative on every element of length between $L+2$ and $k$. Let $$-M= \min_{w \in F_n,\; L+2 \leq \ell(w) \leq k} \eta(C(\gamma_w)).$$ Let $S$ be a finite set of elements of $\bigcup_{i=1}^r A_i$ such that for every element $w \in S_k$, there exists $g_w \in S$ such that $g_w$ is an extension of $w$. The set exists by Lemma~\ref{Lem sets S_k}~$(2)$. Let $$\eta_{\mathcal{A}}=\sum_{g \in S} \eta_{[g]}.$$ By Lemma~\ref{Lem finite set of words determines vertex group system}~$(3)$, for every $w \in F_n$ such that $C(\gamma_w) \in \mathrm{Cyl}(\mathscr{C})$, we have $\eta_{\mathcal{A}}(C(\gamma_w))=0$. Moreover for every $w \in \bigcup_{i=L+2}^k S_i$, Lemma~\ref{Lem sets S_k}~$(2)$ implies that there exists $w' \in S_k$ such that $w'$ is an extension of $w$. In particular, for every $w \in \bigcup_{i=L+2}^k S_i$, we have $\eta_{\mathcal{A}}(C(\gamma_w))>0$. By finiteness of $\bigcup_{i=L+2}^k S_i$, there exists a constant $R>0$ such that for every element $w$ in $\bigcup_{i=L+2}^k S_i$, we have $R\,\eta_{\mathcal{A}}(C(\gamma_w)) \geq M$. 

Then $\eta+ R\,\eta_{\mathcal{A}}$ is nonnegative on words of length between $L+2$ and $k$ and coincides with $\eta_0$ on elements $w \in F_n$ such that $C(\gamma_w) \in \mathrm{Cyl}(\mathscr{C})$. This concludes the proof.
\hfill\qedsymbol

\bigskip

\noindent{\it Proof of Proposition~\ref{Prop density rational currents} } The proof follows \cite[Lemma~3.15]{gupta2017relative} (see also~\cite{Martin95}). Let $\mathscr{C}$ be the set defined above Lemma~\ref{Lem finite set of words determines vertex group system}. Let $\eta_0$ be a relative current and let $k \geq L+2$. Note that every word in $\mathscr{C}$ has length at most equal to $k$. Let $P$ be the constant given by Lemma~\ref{Lem approximation signed measured currents by rational currents}. Note that there exists an element $w'$ in $\mathscr{C}$ such that $\eta_0(C(\gamma_{w'})) >0$. By additivity of $\eta_0$, there exists an element $w_0 \in F_n$ with $\ell(w_0)=k$ and $C(\gamma_{w_0}) \in \mathrm{Cyl}(\mathscr{C})$ and such that $\eta_0(C(\gamma_{w_0})) >0$. Let $R>0$ be such that $R\,\eta_0(C(\gamma_{w_0})) > P$. By Lemma~\ref{Lem k extension}, there exists a signed measured current $\eta$ which is a $k$-extension of $\eta_0$. By Lemma~\ref{Lem approximation signed measured currents by rational currents} applied to $R\,\eta$ and $k'=L+2$, there exists $\alpha_1 \in F_n-\{e\}$ such that for every $w \in F_n-\{e\}$ of length between $L+2$ and $k$, we have $$R\,\eta(C(\gamma_w)) \geq \eta_{[\alpha_1]}(C(\gamma_w)).$$ Suppose first that for every $w \in F_n$ of length between $L+2$ and $k$, we have $$R\,\eta(C(\gamma_w)) \leq \eta_{[\alpha_1]}(C(\gamma_w))+P.$$ Then we stop the process and choose $\alpha_1$. Otherwise, we apply Lemma~\ref{Lem approximation signed measured currents by rational currents} to $R\eta-\eta_{[\alpha_1]}$ and $k'=L+2$. This shows that there exists $\alpha_2 \in F_n-\{e\}$ such that for every $w \in F_n-\{e\}$ of length between $L+2$ and $k$, we have $$R\,\eta(C(\gamma_w))-\eta_{[\alpha_1]}(C(\gamma_w)) \geq \eta_{[\alpha_2]}(C(\gamma_w)).$$ Applying these arguments iteratively (the process stops by Remark~\ref{Rmq lem approximation by rational currents}~$(2)$), we see that there exist $\alpha_1,\ldots,\alpha_p \in F_n-\{e\}$ such that for every element $w \in F_n-\{e\}$ of length between $L+2$ and $k$, we have: 
$$ \sum\limits_{i=1}^p \eta_{[\alpha_i]}(C(\gamma_w)) \leq R\,\eta(C(\gamma_w)) \leq \sum\limits_{i=1}^p \eta_{[\alpha_i]}(C(\gamma_w)) +P.$$

We claim that there exists $i \in \{1,\ldots,p\}$ such that $\alpha_i$ is nonperipheral. Indeed, suppose towards a contradiction that for every $i \in \{1,\ldots,p\}$, the element $\alpha_i$ is peripheral. By Lemma~\ref{Lem finite set of words determines vertex group system}~$(3)$, we have $$\sum\limits_{i=1}^p \eta_{[\alpha_i]}(C(\gamma_{w_0}))=0.$$ This implies that $R\,\eta(C(\gamma_{w_0})) \leq P$. This contradicts the construction of $\eta$. Therefore there exists $i \in \{1,\ldots,p\}$ such that $\alpha_i$ is nonperipheral. Let $S$ be the subset of $\{\alpha_1,\ldots,\alpha_p\}$ containing every nonperipheral element. Then, for every element $w \in F_n$ of length $k$ such that $C(\gamma_w) \in \mathrm{Cyl}(\mathscr{C})$ we have:
$$\left| \eta(C(\gamma_w))-\frac{\sum_{\alpha \in S}\eta_{[\alpha]}(C(\gamma_w))}{R} \right| \leq \frac{P}{R}.$$
For $\alpha \in S$, let $\overline{\eta}_{[\alpha]}$ be the restriction of $\eta_{[\alpha]}$ to the Borel subsets of $\partial^2(F_n,\mathcal{A})$. By construction of $\eta$, for every element $w \in F_n$ of length at most $k$ such that $C(\gamma_w) \in \mathrm{Cyl}(\mathscr{C})$, we have:
$$\left| \eta_0(C(\gamma_w))-\frac{\sum_{\alpha \in S}\overline{\eta}_{[\alpha]}(C(\gamma_w))}{R} \right| \leq \frac{P}{R}.$$

Since $R$ can be chosen arbitrarily large, we can approximate relative currents by sum of rational relative currents. For $m \in \NN^*$, let $\beta^m= \prod_{\alpha \in S} \alpha^m$ (for any total order on $S$). Then there exists $m \in \NN^*$ such that $\sum_{\alpha \in S} \eta_{[\alpha_i]}$ can be approximated by $\frac{1}{m}\eta_{[\beta^m]}$. This concludes the proof.
\hfill\qedsymbol

\bibliographystyle{alphanum}
\bibliography{bibliographie}

\begin{thebibliography}{Gue2}

\bibitem[AM]{AbbottManning2021}
C.~Abbott and J.~Manning.
\newblock {\it Acylindrically hyperbolic groups and their quasi-isometrically
  embedded subgroups}.
\newblock {Preprint {\tt [arXiv:2105.02333]}}.

\bibitem[BH]{BesHan92}
M.~Bestvina and M.~Handel.
\newblock {\it Train tracks and automorphisms of free groups}.
\newblock {Ann. of Math. {\bf 135} (1992) 1--51}.

\bibitem[Bon1]{Bonahon86}
F.~Bonahon.
\newblock {\it Bouts des vari{\'e}t{\'e}s hyperboliques de dimension 3}.
\newblock {Ann. of Math. {\bf 124} (1986) 71--158}.

\bibitem[Bon2]{Bonahon91}
F.~Bonahon.
\newblock {\it Geodesic currents on negatively curved groups}.
\newblock {In ``Arboreal group theory'', R. Alperin ed., pp. 143--168, Pub.
  M.S.R.I. {\bf 19}, Springer Verlag 1991}.

\bibitem[Bou]{Bourbaki65}
N.~Bourbaki.
\newblock {\it Int\'egration : chap.~1, 2, 3 et 4}.
\newblock {Hermann, 1965}.

\bibitem[Bow]{Bowditch2012}
B.~Bowditch.
\newblock {\it Relatively hyperbolic groups}.
\newblock {Internat. J. Algebra Comput. (3) {\bf 22} (2012)}.

\bibitem[CHL]{CouHilLus08III}
T.~Coulbois, A.~Hilion, and M.~Lustig.
\newblock {\it $\mathbb{R}$-trees and laminations for free groups. III.
  Currents and dual $\mathbb{R}$-tree metrics}.
\newblock {J. London Math. Soc. (2) {\bf 78} (2008) 755--766}.

\bibitem[Coh]{Cohn80}
D.~Cohn.
\newblock {\it Measure theory}.
\newblock {Birkh\"auser, 1980}.

\bibitem[DT]{DowdallTaylor2018}
S.~Dowdall and S.J. Taylor.
\newblock {\it Hyperbolic extensions of free groups}.
\newblock {Geom. Topol. {\bf 22} (2018) 517--570}.

\bibitem[EU]{ErlanUya2020}
V.~Erlandsson and C.~Uyanik.
\newblock {\it Length functions on currents and applications to dynamics and
  counting}.
\newblock {In the tradition of Thurston, 423–-458, Springer, Cham, 2020}.

\bibitem[GH]{Guirardelhorbez19laminations}
V.~Guirardel and C.~Horbez.
\newblock {\it Algebraic laminations for free products and arational trees}.
\newblock {Alg. Geom. Topol. (5) {\bf 19} (2019) 2283--2400}.

\bibitem[GL]{GabLev95}
D.~Gaboriau and G.~Levitt.
\newblock {\it The rank of actions on $\RR$-trees}.
\newblock {Ann. Scien. Ec. Norm. Sup. (4) {\bf 28} (1995) 549--570}.

\bibitem[Gue1]{Guerch2021NorthSouth}
Y.~Guerch.
\newblock {\it North-South type dynamics of relative atoroidal automorphisms of
  free groups}.
\newblock {In preparation}.

\bibitem[Gue2]{Guerch2021Polygrowth}
Y.~Guerch.
\newblock {\it Polynomial growth and subgroups of $\mathrm{Out}(F_n)$}.
\newblock {In preparation}.

\bibitem[Gup]{gupta2017relative}
R.~Gupta.
\newblock {\it Relative currents}.
\newblock {Conform. Geom. Dyn. Amer. Math. Soc. {\bf 21} (2017) 319--352}.

\bibitem[HM]{HandelMosher20}
M.~Handel and L.~Mosher.
\newblock {\it Subgroup Decomposition in $\mathrm{Out}(F_n)$}.
\newblock {Mem. Amer. Math. Soc. {\bf 264} (2020)}.

\bibitem[Kap]{Kapovich2006}
I.~Kapovich.
\newblock {\it Currents on free groups}.
\newblock {In Topological and asymptotic aspects of group theory, Contemp.
  Math. {\bf 394} (2006) Amer. Math. Soc., Providence, RI, 149--176}.

\bibitem[KL]{KapLus09}
I.~Kapovich and M.~Lustig.
\newblock {\it Geometric intersection number and analogues of the curve complex
  for free groups}.
\newblock {Geom. \& Topo. {\bf 13} (2009) 1805--1833}.

\bibitem[KR]{KapovichRafi14}
I.~Kapovich and K.~Rafi.
\newblock {\it On hyperbolicity of free splitting and free factor complexes}.
\newblock {Groups Geom. Dyn. (2) {\bf 8} (2014) 391--414}.

\bibitem[Lev]{Levitt09}
G.~Levitt.
\newblock {\it Counting growth types of automorphisms of free groups}.
\newblock {Geom. Funct. Anal. {\bf 19} (2009) 1119--1146}.

\bibitem[Mar]{Martin95}
R.~Martin.
\newblock {\it Non-uniquely ergodic foliations of thin-type, measured currents
  and automorphisms of free groups}.
\newblock {PhD thesis, University of California, Los Angeles, 1995}.

\bibitem[RS]{RuelleSullivan1975}
D.~Ruelle and D.~Sullivan.
\newblock {\it Currents, flows and diffeomorphisms}.
\newblock {Topology {\bf 14} (1975) 319-–327}.

\bibitem[Swe]{Swenson97}
E.L. Swenson.
\newblock {\it Quasi-convex groups of isometries of negatively curved spaces}.
\newblock {Geometric topology and geometric group theory (Milwaukee, WI, 1997).
  Topology Appl. {\bf 110} (2001) 119--129}.

\bibitem[Thu]{Thurston88}
W.~Thurston.
\newblock {\it On the geometry and dynamics of diffeomorphisms of surfaces}.
\newblock {Bull. Amer. Math. Soc. {\bf 19} (1988) 417-432}.

\bibitem[Tra]{Tran19}
H.C. Tran.
\newblock {\it On strongly quasiconvex subgroups}.
\newblock {Geom. Topol. {\bf 23} (2019) 1173--1235}.

\bibitem[Uya1]{Uyanik2015}
C.~Uyanik.
\newblock {\it Generalized north-south dynamics on the space of geodesic
  currents}.
\newblock {Geom. Dedicata {\bf 177} (2015) 129--148}.

\bibitem[Uya2]{Uyanik2019}
C.~Uyanik.
\newblock {\it Hyperbolic extensions of free groups from atoroidal ping-pong}.
\newblock {Algebr. Geom. Topol. (3) {\bf 19} (2019) 1385--1411}.

\end{thebibliography}

\noindent \begin{tabular}{l}
Laboratoire de mathématique d'Orsay\\
UMR 8628 CNRS \\
Université Paris-Saclay\\
91405 ORSAY Cedex, FRANCE\\
{\it e-mail: yassine.guerch@universite-paris-saclay.fr}
\end{tabular}

\end{document}